\documentclass[12pt,reqno]{amsart}

\usepackage{amsmath,amsthm,amssymb,mathtools,enumerate,ytableau,graphicx,mathrsfs}
\usepackage[all]{xy}
\usepackage{aliascnt,mathdots}
\usepackage{xcolor}
\usepackage{hyperref}
\hypersetup{
    colorlinks=true,
    linkcolor=blue,
    citecolor=blue,
    urlcolor=blue,
    bookmarksnumbered=true,
    bookmarksopen=true,
}
\usepackage[margin=1in]{geometry}
\allowdisplaybreaks

\newtheorem{theorem}{Theorem}[section]

\newaliascnt{proposition}{theorem}
\newtheorem{proposition}[proposition]{Proposition}
\aliascntresetthe{proposition}

\newaliascnt{lemma}{theorem}
\newtheorem{lemma}[lemma]{Lemma}
\aliascntresetthe{lemma}

\newaliascnt{corollary}{theorem}
\newtheorem{corollary}[corollary]{Corollary}
\aliascntresetthe{corollary}

\newtheorem{thmA}{Main Theorem}

\numberwithin{equation}{section}

\theoremstyle{definition}
\newaliascnt{definition}{theorem}
\newtheorem{definition}[definition]{Definition}
\aliascntresetthe{definition}

\newaliascnt{example}{theorem}
\newtheorem{example}[example]{Example}
\aliascntresetthe{example}

\newaliascnt{remark}{theorem}
\newtheorem{remark}[remark]{Remark}
\aliascntresetthe{remark}
\newtheorem*{unnumberedremark}{Remark}

\mathtoolsset{showonlyrefs}

\newcommand{\N}{\mathbb{N}}
\newcommand{\Z}{\mathbb{Z}}

\newcommand{\R}{\mathbb{R}}
\newcommand{\sD}{\mathscr{D}}
\newcommand{\sA}{\mathscr{A}}
\newcommand{\bk}{\boldsymbol{k}}
\newcommand{\bm}{\boldsymbol{m}}
\newcommand{\bn}{\boldsymbol{n}}

\newcommand{\bN}{\boldsymbol{N}}
\newcommand{\ba}{\boldsymbol{a}}
\newcommand{\bb}{\boldsymbol{b}}

\newcommand{\bkappa}{\boldsymbol{\kappa}}
\DeclareMathOperator{\sgn}{sgn}
\DeclareMathOperator{\tr}{Tr}
\DeclareMathOperator{\dep}{dep}
\DeclareMathOperator{\wt}{wt}
\DeclareMathOperator{\id}{id}
\DeclareMathOperator{\SSYT}{SSYT}
\DeclareMathOperator{\Cor}{Cor}

\newcommand{\vsmall}{\rotatebox[origin=c]{-90}{$\scriptstyle<$}}

\newlength{\tableauCell}
\newlength{\tableauRule}
\AtBeginDocument{%
  \setlength{\tableauCell}{1.2em}%
  \setlength{\tableauRule}{0.04em}%
}
\ytableausetup{centertableaux,boxsize=\tableauCell,boxframe=\tableauRule}

\newcommand{\tableauparens}[1]{\left(\ #1\ \right)}

\newenvironment{resizedtableau}
  {\begin{ytableau}[\color{black}]}
  {\end{ytableau}}

\newcommand{\tableauSeam}{%
  $\vcenter{\hrule width \tableauRule
    height \dimexpr\tableauCell+\tableauRule\relax}$}
\newcommand{\leftSeam}{\rlap{\kern.8pt\tableauSeam}}
\newcommand{\rightSeam}{\llap{\tableauSeam\kern.8pt}}
\newcommand{\seamCell}[3]{\hbox to \tableauCell{#1\hss$#2$\hss#3}}

\newcommand{\cyctabdiag}[2]{%
  \begin{resizedtableau}
    #1 & \seamCell{}{#2}{\rightSeam}\\
    \seamCell{\leftSeam}{#2}{} & #1
  \end{resizedtableau}}
\newcommand{\cyctabrow}[2]{%
  \begin{resizedtableau}
    \seamCell{\leftSeam}{#1}{} & \seamCell{}{#2}{\rightSeam}\\
    \seamCell{\leftSeam}{#2}{} & \seamCell{}{#1}{\rightSeam}
  \end{resizedtableau}}

\newcommand{\cyctaballtwo}{%
  \vcenter{\hbox{$\scriptstyle n$}}\,{\left\{
  \begin{resizedtableau}
    \seamCell{\leftSeam}{2}{} & \seamCell{}{2}{\rightSeam}\\
    \seamCell{\leftSeam}{\vcenter{\hbox{\scalebox{.65}{$\vdots$}}}}{} &
      \seamCell{}{\vcenter{\hbox{\scalebox{.65}{$\vdots$}}}}{\rightSeam}\\
    \seamCell{\leftSeam}{2}{} & \seamCell{}{2}{\rightSeam}
  \end{resizedtableau}\right.}}

\newcommand{\cycsingle}{%
  \begin{resizedtableau}\seamCell{\leftSeam}{1}{\rightSeam}\end{resizedtableau}}
\newcommand{\cycgallery}{%
  \begin{resizedtableau}
    \none & 1 & \seamCell{}{2}{\rightSeam}\\
    \seamCell{\leftSeam}{2}{} & 1 & \seamCell{}{1}{\rightSeam}\\
    \seamCell{\leftSeam}{2}{} & 2
  \end{resizedtableau}}
\newcommand{\cycsix}{%
  \begin{resizedtableau}
    2 & \seamCell{}{2}{\rightSeam}\\
    \seamCell{\leftSeam}{2}{} & \seamCell{}{2}{\rightSeam}\\
    \seamCell{\leftSeam}{2}{} & 2
  \end{resizedtableau}}

\newcommand{\cyctworow}{%
  \begin{resizedtableau}
    \seamCell{\leftSeam}{2}{} & \seamCell{}{2}{\rightSeam}
  \end{resizedtableau}}

\newcommand{\periodicEntry}[1]{\resizebox{.9\tableauCell}{!}{$m_{#1}$}}

\title[A cyclic analogue of Schur multiple zeta values]{A cyclic analogue of Schur multiple zeta values}

\makeatletter
\let\original@setaddresses\@setaddresses
\renewcommand{\@setaddresses}{%
  \begingroup
  \raggedright
  \parindent=0pt
  \original@setaddresses
  \endgroup
}
\makeatother

\author{Ku-Yu Fan}
\address{Graduate School of Mathematics, Nagoya University, Nagoya, Japan.}
\email{ku-yu.fan.d2@math.nagoya-u.ac.jp}

\author{Hau-Yi Lin}
\address{Department of Applied Mathematics, National Yang Ming Chiao Tung University, Hsinchu~30010, Taiwan.}
\email{zaq1bgt5cde3mju7@gmail.com}

\author{Jinbo Yu}
\address{Graduate School of Mathematics, Nagoya University, Nagoya, Japan.}
\email{jinbo.yu.e6@math.nagoya-u.ac.jp}


\usepackage[nameinlink]{cleveref}
\crefname{theorem}{Theorem}{Theorems}
\crefname{proposition}{Proposition}{Propositions}
\crefname{lemma}{Lemma}{Lemmas}
\crefname{corollary}{Corollary}{Corollaries}
\crefname{definition}{Definition}{Definitions}
\crefname{example}{Example}{Examples}
\crefname{remark}{Remark}{Remarks}
\crefname{thmA}{Main Theorem}{Main Theorems}
\crefname{section}{Section}{Sections}
\Crefname{section}{Section}{Sections}
\crefname{subsection}{Section}{Sections}
\Crefname{subsection}{Section}{Sections}
\crefname{equation}{}{}
\creflabelformat{equation}{(#2#1#3)}
\begin{document}
\date{\today}

\begin{abstract}
In this paper we introduce a cyclic analogue of Schur multiple zeta values. For a diagonally constant index on a periodically extended skew Young diagram, it is a sum over semi-standard cyclic Young tableaux. It includes the cyclic multiple zeta values of Hirose, Murahara and Murakami as a special case. The first two main results are a cyclic analogue of the Maesaka--Seki--Watanabe formula for truncated sums and an {integral--series identity} under an admissibility condition. As an application, we give a proof of $4^n\zeta(\{1,3\}^n)=\zeta(\{4\}^n)$ by counting linear extensions of a poset. The third main result expresses a cyclic Schur multiple harmonic sum as a finite signed sum of Schur multiple harmonic sums with the same number of boxes and the same multiset of exponents.
\end{abstract}

\maketitle

\begingroup
  \tableofcontents
\endgroup

\section{Introduction}\label{sec:intro}

The purpose of this paper is to introduce a cyclic analogue of Schur multiple zeta values and to {prove its integral--series identity} via its Maesaka--Seki--Watanabe formula. For $r\geq1$ and $(k_1,\ldots,k_r)\in\N^r$ with $k_r\geq2$, the multiple zeta values and the multiple zeta-star values are defined by {the following absolutely convergent} series:
\begin{equation}\label{eq:defMZV}
\zeta(k_1,\ldots,k_r)=\sum_{0<n_1<\cdots<n_r}\frac{1}{n_1^{k_1}\cdots n_r^{k_r}},\qquad
\zeta^\star(k_1,\ldots,k_r)=\sum_{0<n_1\leq\cdots\leq n_r}\frac{1}{n_1^{k_1}\cdots n_r^{k_r}}.
\end{equation}
Multiple zeta values {also admit iterated integral expressions} on $[0,1]$ \cite{LM,Za} (see also \cref{thm:iterated-integral}). {These expressions} give the shuffle relation of multiple zeta values \cite{Ka}.

Schur multiple zeta values were introduced by Nakasuji, Phuksuwan and Yamasaki \cite{NPY}. They generalize multiple zeta values by replacing the index $(k_1,\ldots,k_r)$ by {an index on a Young diagram}. For example, for integers $a\geq2$ and $b\geq1$,
\[
\zeta\tableauparens{\begin{ytableau}a&b\\b&a\end{ytableau}}=\sum_{\begin{smallmatrix}
\mathllap{0<{}}m_{11}&\leq&m_{12}\\
\vsmall&&\vsmall\\
m_{21}&\leq&m_{22}
\end{smallmatrix}}\frac{1}{m_{11}^{a}m_{12}^{b}m_{21}^{b}m_{22}^{a}}.
\]
A single column gives a multiple zeta value and a single row gives a multiple zeta-star value. The series converges if the exponents at the corners are at least $2$ \cite{NPY}. An index is called \emph{diagonally constant} if its exponents are constant along diagonals. For such indices, Hirose, Murahara and Onozuka \cite[Theorem~1.2]{HMO} expressed Schur multiple zeta values as Yamamoto's poset iterated integrals \cite{Y1} (see also \cref{thm:hmo}).

For $r\geq1$, an index $\bk=(k_1,\ldots,k_r)\in\N^r$ and a positive integer $N$, the multiple harmonic sum is
\begin{equation}\label{eq:intro-mhs}
\zeta_{<N}(\bk)=\sum_{0<n_1<\cdots<n_r<N}\frac{1}{n_1^{k_1}\cdots n_r^{k_r}}.
\end{equation}
The Maesaka--Seki--Watanabe formula \cite[Theorem~1.3]{MSW} (see also \cref{thm:msw}) states that
\begin{equation}\label{eq:intro-msw}
\zeta_{<N}(\bk)=\zeta^\flat_{<N}(\bk)\qquad(N\geq1),
\end{equation}
where $\zeta^\flat_{<N}(\bk)$ is the $\flat$-sum recalled in \cref{sec:msw}. For example,
\[
\zeta^\flat_{<N}(2,3)=\sum_{0<n_1\leq n_2<n_3\leq n_4\leq n_5<N}\frac{1}{(N-n_1)n_2(N-n_3)n_4n_5}.
\]
For $k_r\geq2$, the $\flat$-sum is a modified Riemann sum with step $1/N$ {for the iterated integral in \cref{thm:iterated-integral}}. {Letting $N\to\infty$ in \eqref{eq:intro-msw} identifies this integral with the series for $\zeta(\bk)$.} Yamamoto \cite[Theorem~3.6]{Y2} extended \eqref{eq:intro-msw} to diagonally constant Schur multiple zeta values (see also \cref{thm:y2}). His formula is a discretization of the {iterated} integral expression of Hirose, Murahara and Onozuka \cite[Theorem~1.2]{HMO}.

The cyclic multiple zeta values of Hirose, Murahara and Murakami \cite{HMM} are defined by strictly increasing blocks of summation variables joined by cyclic weak inequalities. These values admit an integral expression over a cyclic domain \cite[Theorem~1]{HMM} (see also \cref{thm:hmm}).

Let $\lambda=(\lambda_1,\ldots,\lambda_r)$ be a partition and $\mu\subseteq\lambda$. For $q\geq0$ we extend $\lambda/\mu$ periodically by $(i,j)\longmapsto(i+q,j-\lambda_1)$. If the extension is again a skew diagram, we call the quotient by this translation a cyclic Young diagram and denote it by $(\lambda/\mu)_q$ (\cref{def:cyclic-diagram}). A semi-standard cyclic Young tableau is a periodic filling which is weakly increasing along rows and strictly increasing down columns. The diagonal classes of $(\lambda/\mu)_q$ are indexed by the residue of $i-j$ modulo $s=\lambda_1+q$. Let $\bk=(k_{i,j})_{(i,j)\in\lambda/\mu}$ be a positive integer index which is constant along the diagonal classes. Such an index is again called diagonally constant. Its cyclic Schur multiple harmonic sum $\zeta_{<N}(\bk)$ is the corresponding weighted sum over semi-standard cyclic Young tableaux with entries less than $N$.

An index on a cyclic Young diagram is drawn as a filled diagram, and the double edges show how the diagram is glued. For example, for integers $a,b\geq1$ and $q=1$ we have
\begin{equation}\label{eq:intro-2x2}
\zeta_{<N}\tableauparens{\cyctabdiag{a}{b}}=\sum_{\substack{\begin{smallmatrix}
\mathllap{0<{}}m_{11}&\leq&m_{12}\\
\vsmall&&\vsmall\\
m_{21}&\leq&m_{22}\mathrlap{{}<N,}
\end{smallmatrix}\\
\begin{smallmatrix}m_{12}&\leq&m_{21}\end{smallmatrix}}}\frac{1}{m_{11}^{a}m_{12}^{b}m_{21}^{b}m_{22}^{a}}.
\end{equation}
This is the sum over the semi-standard tableaux of shape $(2,2)$ with the additional relation $m_{12}\leq m_{21}$ coming from the periodicity. For $q=0$ each row is closed into a cycle, and
\[
\zeta_{<N}\tableauparens{\cyctabrow{a}{b}}=\sum_{\substack{\begin{smallmatrix}
\mathllap{0<{}}m_{11}&\leq&m_{12}\\
\vsmall&&\vsmall\\
m_{21}&\leq&m_{22}\mathrlap{{}<N,}
\end{smallmatrix}\\
\begin{smallmatrix}m_{12}&\leq&m_{11},&m_{22}&\leq&m_{21}\end{smallmatrix}}}\frac{1}{m_{11}^{a}m_{12}^{b}m_{21}^{b}m_{22}^{a}}
=\zeta_{<N}(a+b,a+b).
\]
Here the conditions force $m_{11}=m_{12}$ and $m_{21}=m_{22}$. The first main result is a cyclic analogue of \eqref{eq:intro-msw}.

\begin{thmA}\label{thmA:finite}
For every diagonally constant index $\bk$ on a cyclic Young diagram and every positive integer $N$, we have
\[
\zeta_{<N}(\bk)=\zeta^\flat_{<N}(\bk).
\]
\end{thmA}
Here $\zeta^\flat_{<N}(\bk)$ is the cyclic Schur $\flat$-sum of \cref{def:cyclic-flat}. It is obtained by expanding every diagonal class of exponent $k$ into $k$ layers, as in Yamamoto's $\flat$-sum \cite[Definition~3.5]{Y2} (see also \eqref{eq:schur-flat}). \Cref{thmA:finite} is \cref{thm:finite}, proved in \cref{sec:transfer}. In diagonal coordinates, both finite sums are traces of products of transfer matrices (see \eqref{eq:trace-series} and \eqref{eq:trace-flat}). The transport identities in \cref{lem:transport} show that the two products are {similar}. They follow from lemmas of Yamamoto \cite[Lemmas~3.11 and~3.12]{Y2} (see also \cref{lem:yamamoto-transport}).

The second main result is an {integral--series identity for cyclic Schur multiple zeta values}. For the limit $N\to\infty$ we use an admissibility condition. A \emph{cyclic corner} is the class of a box with neither a right nor a lower neighbor in the periodic diagram. If cyclic corners exist, we require the exponent at every cyclic corner to be at least $2$. If there is no cyclic corner, then $q=0$, $\mu=\varnothing$ and $\lambda=(\lambda_1^r)$ by \cref{lem:no-corner}. In this case the row sum $K=\sum_{j=1}^{\lambda_1}k_{i,j}$ is independent of $i$, and we require $K\geq2$. An index with these properties is called \emph{admissible}. In \cref{sec:integral} we associate an admissible $2$-labeled poset $X_{\bk}$ to every admissible cyclic index. Let $I(X_{\bk})$ be its Yamamoto integral \cite[Definition~2.1]{Y1} (see also \eqref{eq:yamamoto-integral}).

\begin{thmA}\label{thmA:integral}
Let $\bk$ be an admissible diagonally constant index on a cyclic Young diagram. Then the series $\zeta(\bk)$ converges{, and we have the integral--series identity}
\[
\zeta(\bk)=I(X_{\bk}).
\]
\end{thmA}
\Cref{thmA:integral} is \cref{thm:integral} in \cref{sec:integral}. The passage to the limit uses \cref{lem:approx}. The proof of \cref{lem:approx} applies the estimates of Seki \cite[Lemma~2.2 and Proposition~2.3]{Se} (see also \cref{lem:seki-estimate,lem:seki-strict}) to the contributions with equal summation variables. For example, for $\cyctabdiag{a}{b}$ the condition of \cref{thmA:integral} is $a\geq2$. Cyclic multiple zeta values occur when each diagonal class consists of a single box (\cref{rem:cyclic-hmm}). In this case \cref{thmA:integral} recovers the {cyclic integral--series identity} of Hirose, Murahara and Murakami \cite[Theorem~1]{HMM} (see also \cref{thm:hmm}).

As an application, we prove in \cref{cor:zagier} that
\[
4^n\zeta(\{1,3\}^n)=\zeta(\{4\}^n)
\]
for every integer $n\geq0$. Here $\{1,3\}^n$ denotes $n$ repetitions of $1,3$, and $\zeta(\emptyset)=1$. This identity was conjectured by Zagier \cite{Za} and first proved by Borwein, Bradley, Broadhurst and Lison\v{e}k \cite{BBBL01}. They also gave a combinatorial proof using the shuffle product \cite[Theorem~1]{BBBL98}. Our proof uses the cyclic $n\times2$ rectangle $((2^n)/\varnothing)_0$ with exponent $2$ on every box. Its cyclic Young tableau conditions identify the two entries in each row, so its truncated sum is $\zeta_{<N}(\{4\}^n)$. Its associated poset has exactly $4^n$ linear extensions, all with the label word $(\bullet\bullet\circ\circ)^n$. Hence its integral is $4^n\zeta(\{1,3\}^n)$ by Yamamoto's formula \cite[Corollary~2.4]{Y1} (see also \cref{lem:linear-extension}).

In \cref{sec:finite-mzv} we apply \cref{thmA:finite} with $N$ a prime ${p}$ and reduce modulo ${p}$. The resulting \cref{cor:finite-A} expresses the finite analogues of cyclic Schur multiple zeta values as $\Z$-linear combinations of finite multiple zeta values \cite[Definition~7.1]{Ka}. These finite analogues include Murahara's cyclic analogues of finite multiple zeta values \cite[Definition~1.10]{Mur}.

The third main result expresses a cyclic Schur multiple harmonic sum {as a signed sum of} Schur multiple harmonic sums. For $q>0$, the cyclic Young diagram $(\lambda/\mu)_q$ is a cylindric diagram of type $(q,s)$ in the sense of Postnikov \cite{Po}.

\begin{thmA}\label{thmA:decomp}
Let $\bk$ be a diagonally constant index on $(\lambda/\mu)_q$ with no empty diagonal class. Then for every positive integer $N$,
\[
\zeta_{<N}(\bk)=\sum_{\nu\in\mathcal A}\sgn(\pi_\nu)\,\zeta_{<N}(\nu/\mu;\bkappa),
\]
where the finite set $\mathcal A$ of partitions, the permutations $\pi_\nu$ and the periodic sequence $\bkappa$ of exponents are given in \cref{sec:decomp}. Every summand has the same number of boxes and the same multiset of exponents as $\bk$.
\end{thmA}
The full statement is \cref{thm:decomp} in \cref{sec:decomp}. For $q>0$, the same set of semi-standard cyclic Young tableaux with the weight $x^T=\prod_{(i,j)}x_{T(i,j)}$ gives a cylindric skew Schur polynomial in $x_1,\ldots,x_{N-1}$. Its expansion into skew Schur polynomials goes back to the lattice path formula of Gessel and Krattenthaler \cite[Proposition~1]{GK} (see also \cite[Theorems~6.3 and~6.5]{McN}). Our proof is direct. A tableau corresponds to a family of lattice paths, and the proof switches the tails of two paths at their last meeting modulo $s$. Gessel and Krattenthaler prove their formula by a similar switch. For example, for integers $a,b\geq1$,
\[
\zeta_{<N}\tableauparens{\cyctabdiag{a}{b}}
=\zeta_{<N}\tableauparens{\begin{ytableau}a&b\\b&a\end{ytableau}}-\zeta_{<N}\tableauparens{\begin{ytableau}a\\b\\b\\a\end{ytableau}}.
\]
The same decomposition holds for the $\flat$-sums (\cref{rem:flat-decomp}).

In \cref{sec:app-motivic} we use the first named author's motivic version \cite{Fa} of Yamamoto's integrals in \cref{thmA:integral} to define motivic cyclic Schur multiple zeta values. If all the indices involved are admissible, {\cref{cor:motivic-decomp} gives a motivic version of \cref{thmA:decomp}.}

The paper is organized as follows. In \cref{sec:background,sec:msw} we recall multiple, Schur multiple and cyclic multiple zeta values, Yamamoto's integral and the finite identities of Maesaka--Seki--Watanabe and Yamamoto. In \cref{sec:cyclic} we define cyclic Young diagrams and tableaux and prove \cref{thmA:finite}. \Cref{sec:integral,sec:decomp} contain the statements and proofs of \cref{thmA:integral,thmA:decomp}, respectively. The applications in \cref{sec:applications} are finite multiple zeta values, a proof of Zagier's identity and motivic cyclic Schur multiple zeta values.

\vspace{6pt}
{\bf Acknowledgment.} The authors would like to thank Shin-ichiro Seki for his talk at Nagoya University. He introduced us to the MSW formula, which provided the impetus for the main idea of this paper. The authors would like to also thank Shuji Yamamoto for his talk at Kindai University.

\vspace{6pt}
{\bf AI \& computational resource disclosure.} This project began in late 2024, when we noticed the \cref{thmA:integral}. However, for various reasons, it wasn't fully completed until this year. These results, their formulations, and the proof strategies were all developed by human authors. ChatGPT~4 and 5 were used throughout as research assistants: they were used to write Python and Mathematica code to help calculate and verify explicit combinatorial examples. This helped us refine some of the details of the proof. ChatGPT~6 Sol and Claude Opus~5.5 helped to improve the English and draw the pictures in the paper.

\section{Multiple zeta values}\label{sec:background}

In this section we recall multiple zeta values, Schur multiple zeta values and cyclic multiple zeta values, together with their integral expressions.

An \emph{index} is a tuple $\bk=(k_1,\ldots,k_r)\in\N^r$ with $r\geq0$. Its \emph{weight} is $\wt(\bk)=k_1+\cdots+k_r$ and its \emph{depth} is $\dep(\bk)=r$. The index of depth $0$ is the \emph{empty index} $\emptyset$. An index is \emph{admissible} if it is empty or if $r\geq1$ and $k_r\geq2$. The value of every sum and integral defined below at the empty index, the empty diagram or the empty poset is $1$. In particular, $\zeta(\emptyset)=\zeta^\star(\emptyset)=1$.

For a nonempty index $\bk$ of depth $r$, the condition $k_r\geq2$ implies
\[
0\leq\zeta(\bk)\leq\zeta^\star(\bk)
\leq\sum_{n=1}^{\infty}\frac{(1+\log n)^{r-1}}{n^2}<\infty.
\]
Indeed, $\sum_{m=1}^n m^{-1}\leq1+\log n$. If $k_r=1$, fixing $n_1<\cdots<n_{r-1}$ leaves a divergent harmonic tail. Thus convergence is equivalent to $k_r\geq2$.

\begin{theorem}[{\cite{LM}}, {\cite[Theorem~0.1]{MSW}}]\label{thm:iterated-integral}
Let $\bk=(k_1,\ldots,k_r)\in\N^r$ with $r\geq1$ and $k_r\geq2$, and let $w=k_1+\cdots+k_r$. Set $\omega_0(t)=dt/t$, $\omega_1(t)=dt/(1-t)$ and
\[
a_i=\begin{cases}1,&i=1+\sum_{j=1}^{h}k_j\text{ for some }0\leq h<r,\\0,&\text{otherwise}.
\end{cases}
\]
Then the integral converges and
\begin{equation}\label{eq:iterated-integral}
\zeta(k_1,\ldots,k_r)=\int_{0<t_1<\cdots<t_w<1}\prod_{i=1}^w\omega_{a_i}(t_i).
\end{equation}
\end{theorem}
The formula \eqref{eq:iterated-integral} {gives the classical} iterated integral expression {for} multiple zeta values. L\^e and Murakami \cite{LM} recorded this formula, attributed to Kontsevich and Drinfeld. The integrand of \eqref{eq:iterated-integral} is encoded by the word $\bullet\circ^{k_1-1}\cdots\bullet\circ^{k_r-1}$, where $\bullet=\omega_1$ and $\circ=\omega_0$.

\begin{example}\label{ex:iterated-word}
Take $\bk=(1,2)$, of weight $3$. In \cref{thm:iterated-integral}, the labels are $a_1=a_2=1$ and $a_3=0$. Therefore its word is $\bullet\bullet\circ$, and the formula reads
\[
\zeta(1,2)=\int_{0<t_1<t_2<t_3<1}
\underbrace{\frac{dt_1}{1-t_1}}_{\bullet}\,
\underbrace{\frac{dt_2}{1-t_2}}_{\bullet}\,
\underbrace{\frac{dt_3}{t_3}}_{\circ}.
\]
\end{example}

Yamamoto's integral \cite[Definition~2.1]{Y1} replaces the total order in \eqref{eq:iterated-integral} by a partial order. A \emph{$2$-labeled poset} is a pair $(X,\delta_X)$ of a finite partially ordered set $X$ and a map $\delta_X\colon X\to\{\circ,\bullet\}$. It is called \emph{admissible} if $\delta_X({x})=\circ$ for every maximal element ${x}$ of $X$ and $\delta_X({x'})=\bullet$ for every minimal element ${x'}$ of $X$. For an admissible $2$-labeled poset $X$, \emph{Yamamoto's integral} is
\begin{equation}\label{eq:yamamoto-integral}
I(X)=\int_{\Delta_X}\prod_{{x}\in X}\omega_{\delta_X({x})}(t_{x}),\qquad
\Delta_X=\{(t_{x})_{{x}\in X}\in(0,1)^X\mid t_{x}<t_{x'}\text{ if }{x}<{x'}\},
\end{equation}
where $\omega_\circ=\omega_0$ and $\omega_\bullet=\omega_1$.

A \emph{linear extension} of $X$ is a list ${x}_1,{x_2,}\ldots,{x}_w$ of all elements of $X$ such that ${x}_i<{x}_j$ in $X$ implies $i<j$.

\begin{lemma}[{\cite[Corollary~2.4]{Y1}}]\label{lem:linear-extension}
Let $X$ be a nonempty admissible $2$-labeled poset. For a linear extension ${x}_1,{x_2,}\ldots,{x}_w$ of $X$, the word $\delta_X({x}_1)\cdots\delta_X({x}_w)$ is of the form $\bullet\circ^{k_1-1}\cdots\bullet\circ^{k_r-1}$ with $k_r\geq2$. Then $I(X)$ is the sum of $\zeta(k_1,\ldots,k_r)$ over all linear extensions of $X$. In particular, $I(X)$ converges.
\end{lemma}

A \emph{partition} of a natural number $n$ is a tuple $\lambda=(\lambda_1,\ldots,\lambda_l)$ of positive integers satisfying
\[
\lambda_1\geq\cdots\geq\lambda_l\geq1,\qquad|\lambda|:=\lambda_1+\cdots+\lambda_l=n.
\]
Here we also allow the empty partition $\varnothing$, with $|\varnothing|=0$. When comparing two partitions, we extend them by adding zero parts at the end if necessary. The \emph{conjugate} of $\lambda$ is the partition $\lambda'$ with $\lambda'_j=\#\{i\mid\lambda_i\geq j\}$ for $j\geq1$.

Let $\lambda=(\lambda_1,\ldots,\lambda_l)$ and $\mu$ be partitions. Write $\mu\subseteq\lambda$ if $0\leq\mu_i\leq\lambda_i$ for all $i$, after adding zero parts to $\mu$ if necessary. For $\mu\subseteq\lambda$, the \emph{skew Young diagram} of shape $\lambda/\mu$ is the set
\[
\{(i,j)\in\Z^2\mid1\leq i\leq l,\ \mu_i<j\leq\lambda_i\},
\]
which we also denote by $\lambda/\mu$. Here $(i,j)$ is the box in row $i$ and column $j$. For $\mu=\varnothing$, it is the \emph{Young diagram} of $\lambda$, denoted by $\lambda$. For a skew Young diagram $\lambda/\mu$, let $\SSYT(\lambda/\mu)$ be the set of all \emph{semi-standard Young tableaux} of shape $\lambda/\mu$:
\[
\SSYT(\lambda/\mu)=\{T\colon\lambda/\mu\to\N\mid T(i,j)\leq T(i,j+1),\ T(i,j)<T(i+1,j)\}.
\]
Each inequality is imposed whenever both boxes lie in $\lambda/\mu$. For a positive integer $N$, let $\SSYT_{<N}(\lambda/\mu)$ be the subset of those with all entries less than $N$. A box $(i,j)\in\lambda/\mu$ is called a \emph{corner} if $(i+1,j)\notin\lambda/\mu$ and $(i,j+1)\notin\lambda/\mu$.

A map $\bk\colon\lambda/\mu\to\N$ is called an \emph{index} on $\lambda/\mu$. Write $k_{i,j}=\bk(i,j)$ and draw $\bk$ as a filled diagram. For example,
\[
\bk=\begin{ytableau}\none&1&3\\4&2\end{ytableau}
\]
is the index on $(3,2)/(1)$ with $k_{1,2}=1$, $k_{1,3}=3$, $k_{2,1}=4$ and $k_{2,2}=2$. An index $\bk$ on $\lambda/\mu$ is called \emph{admissible} if $k_{i,j}\geq2$ for every corner $(i,j)$ of $\lambda/\mu$. By Nakasuji, Phuksuwan and Yamasaki \cite{NPY}, the \emph{Schur multiple zeta value} of $\bk$
\begin{equation}\label{eq:defSMZV}
\zeta(\bk)=\sum_{T\in\SSYT(\lambda/\mu)}\prod_{(i,j)\in\lambda/\mu}\frac{1}{T(i,j)^{k_{i,j}}}
\end{equation}
converges absolutely for every admissible index.

An index $\bk$ on $\lambda/\mu$ is called \emph{diagonally constant} if $k_{i,j}$ depends only on $i-j$. For ${g}\in\Z$, we write
\begin{equation}\label{eq:diagonal}
D_{g}=\{(i,j)\in\lambda/\mu\mid i-j={g}\},\qquad
J_{g}=\{j\in\Z\mid (j+{g},j)\in D_{g}\}
\end{equation}
for the ${g}$-th diagonal of $\lambda/\mu$ and for the set of its column indices. Each $J_{g}$ is an interval in $\Z$, possibly empty. If $\lambda/\mu$ is nonempty, let ${g}_0$ and ${g}_1$ be the smallest and the largest ${g}$ with $D_{g}\neq\emptyset$. For a diagonally constant index $\bk$ and ${g}_0\leq {g}\leq {g}_1$, we write $k_{g}=k_{i,j}$ for $(i,j)\in D_{g}$, and we put $k_{g}=1$ if $D_{g}=\emptyset$. A skew Young diagram is called \emph{connected} if its edge-adjacency graph is connected.

\begin{theorem}[{\cite[Theorem~1.2]{HMO}}]\label{thm:hmo}
For an admissible diagonally constant index on a connected skew Young diagram, we have
\begingroup
\small
\begin{align*}
  \zeta\tableauparens{%
  \begin{resizedtableau}
  \none & \none & \none & \none &  k_1   \\
  \none &  k_i & \none & \none[\color{black}\iddots] & \none \\
  \none & \none & \none[\color{black}\ddots] \\
  \none & \none[\color{black}\iddots] & \none & k_i  \\
  k_r & \none & \none
  \end{resizedtableau}%
  }=I\left(\
\begin{xy}
{(0,-6) \ar @{{*}-o} (4,-2)}, 
{(4,-2) \ar @{.o} (8,2)}, 
{(0,-5) \ar @/^2mm/ @{-}^{k_1} (7,2)}, 
{(8,2) \ar @{.} (12,6)},  
{(8,2) \ar @{.} (12,-2)}
\end{xy}
\cdots
\begin{xy}
{(32,4) \ar @{{*}-o} (36,8)},   
{(32,4) \ar @{{*}-o} (36,0)},   
{(32,-4) \ar @{{*}-} (36,0)},   
{(32,-4) \ar @{.} (36,-8)},   
{(32,-12) \ar @{{*}-o} (36,-8)},   
{(32,4) \ar @/_2mm/ @{-}_{l_i} (32,-11)}, 
{(36,8) \ar @{.o} (40,12)},   
{(36,8) \ar @{.o} (40,4)},   
{(36,0) \ar @{.} (40,4)},   
{(36,0) \ar @{.} (40,-4)},   
{(36,-8) \ar @{.o} (40,-4)}, 
{(32,5) \ar @/^2mm/ @{-}^{k_i} (39,12)}, 
\end{xy}\ 
\cdots
\begin{xy}
{(44,2) \ar @{.} (48,-2)},   
{(44,-6) \ar @{.{*}} (48,-2)},   
{(48,-2) \ar @{-o} (52,2)}, 
{(52,2) \ar @{.o} (56,6)}, 
{(48,-1) \ar @/^2mm/ @{-}^{k_r} (55,6)}, 
\end{xy}\ 
\right),
\end{align*}
\endgroup
where $l_i$ is the number of boxes in the $i$-th diagonal. The diagonals are numbered from the top right. Thus the $i$-th diagonal is $D_{{g}_0+i-1}$, $k_i$ stands for $k_{{g}_0+i-1}$ and $r={g}_1-{g}_0+1$.
\end{theorem}
For disconnected shapes, take the disjoint union of the posets of the connected components.

\begin{example}\label{ex:yamamoto}
Hasse diagrams of $2$-labeled posets are drawn with the order increasing from bottom to top, as in \cref{thm:hmo}. The single-row case of \cref{thm:hmo} gives
\begin{align*}
\zeta^\star(1,2)
&=\int_{\substack{0<t_1<t_2<1\\0<t_3<t_2}}
\frac{dt_1}{1-t_1}\frac{dt_2}{t_2}\frac{dt_3}{1-t_3}
=I\left(
\begin{xy}
{(-4,-2) \ar @{{*}-o} (0,2)},
{(4,-2) \ar @{{*}-o} (0,2)}
\end{xy}
\right).
\end{align*}
The two orders $t_1<t_3<t_2$ and $t_3<t_1<t_2$ have the same label word $\bullet\bullet\circ$. Hence \cref{lem:linear-extension} gives $\zeta^\star(1,2)=2\zeta(1,2)$.
\end{example}

Next we recall the cyclic multiple zeta values of Hirose, Murahara and Murakami \cite{HMM}. A \emph{multi-index} is an ordered list
\[
\bk=[(k_{1,1},\ldots,k_{1,r_1}),\ldots,(k_{s,1},\ldots,k_{s,r_s})]
\]
of $s\geq1$ nonempty indices. It is called \emph{admissible} if each index $(k_{i,1},\ldots,k_{i,r_i})$ is admissible or equal to $(1)$, and at least one of them is different from $(1)$. For an admissible multi-index, the \emph{cyclic multiple zeta value} is
\begin{equation}\label{eq:defcmzv}
\zeta(\bk)=\sum_{\substack{n_{i,j}>0\\
n_{1,1}<\cdots<n_{1,r_1}\geq n_{2,1}<\cdots<n_{2,r_2}\geq\cdots\geq n_{s,1}<\cdots<n_{s,r_s}\geq n_{1,1}}}
\prod_{i=1}^s\prod_{j=1}^{r_i}\frac{1}{n_{i,j}^{k_{i,j}}}.
\end{equation}

\begin{theorem}[{\cite[Theorem~1]{HMM}}]\label{thm:hmm}
Let $\bk$ be an admissible multi-index and put $K_i=k_{i,1}+\cdots+k_{i,r_i}$. Then
\[
\zeta(\bk)=\int_{\Delta_{\bk}}\prod_{i=1}^s\prod_{j=1}^{K_i}a_{i,j}(t_{i,j})\,dt_{i,j},
\]
where $a_{i,j}(t)=1/(1-t)$ for $j\in\{1,\,k_{i,1}+1,\,\ldots,\,k_{i,1}+\cdots+k_{i,r_i-1}+1\}$ and $a_{i,j}(t)=1/t$ otherwise. Here $\Delta_{\bk}\subset(0,1)^{K_1+\cdots+K_s}$ is the set of all $(t_{i,j})$ with
\[
t_{1,1}<\cdots<t_{1,K_1}>t_{2,1}<\cdots<t_{2,K_2}>\cdots>t_{s,1}<\cdots<t_{s,K_s}>t_{1,1}.
\]
\end{theorem}
{This is} the cyclic integral--series identity.

\begin{example}\label{ex:cmzv}
For $\bk=[(2),(1)]$, the conditions are $n_{1,1}\geq n_{2,1}\geq n_{1,1}$. Hence
\[
\zeta([(2),(1)])=\sum_{m\geq1}\frac{1}{m^3}=\zeta(3).
\]
The integral in \cref{thm:hmm} is over $t_{1,1}<t_{2,1}<t_{1,2}$ with labels $\bullet,\bullet,\circ$. Thus \cref{thm:hmm,thm:iterated-integral} give $\zeta(3)=\zeta(1,2)$.
\end{example}
\section{The Maesaka--Seki--Watanabe formula}\label{sec:msw}

In this section we recall the finite identity of Maesaka, Seki and Watanabe \cite{MSW} and its extension by Yamamoto \cite{Y2} to diagonally constant Schur indices.

Let $N$ be a positive integer and let $\bk=(k_1,\ldots,k_r)$ be an index. Its \emph{multiple harmonic sum} $\zeta_{<N}(\bk)$ is defined in \eqref{eq:intro-mhs}. Suppose now that $\bk$ is nonempty. Let $w=\wt(\bk)$ and $K_i=k_1+\cdots+k_i$ for $0\leq i\leq r$, where $K_0=0$. Let $S_N(\bk)$ be the set of all $(n_1,\ldots,n_w)\in\Z^w$ such that, for $1\leq i\leq w+1$,
\[
n_{i-1}<n_i\quad\text{if }i-1\in\{K_0,K_1,\ldots,K_r\},\qquad
n_{i-1}\leq n_i\quad\text{otherwise}.
\]
Here $n_0=0$ and $n_{w+1}=N$. The \emph{$\flat$-sum} is
\[
\zeta^\flat_{<N}(\bk)=\sum_{(n_1,\ldots,n_w)\in S_N(\bk)}\prod_{i=1}^r\frac{1}{(N-n_{K_{i-1}+1})\displaystyle\prod_{j=K_{i-1}+2}^{K_i}n_j}.
\]
For $k_i=1$ the product over $j$ is empty. For an admissible index, this sum is a modified Riemann sum for \cref{eq:iterated-integral}.

\begin{theorem}[{\cite[Theorem~1.3]{MSW}}]\label{thm:msw}
For every index $\bk$ and positive integer $N$, we have
\[
\zeta_{<N}(\bk)=\zeta^\flat_{<N}(\bk).
\]
\end{theorem}

\begin{example}
Take $\bk=(2)$ and $N=3$. On the series side,
\[
\zeta_{<3}(2)=1+\frac1{2^2}=\frac54.
\]
The $\flat$-sum uses $1\leq n_1\leq n_2<3$ and gives
\[
\zeta^\flat_{<3}(2)
=\frac1{(3-1)\cdot1}+\frac1{(3-1)\cdot2}
 +\frac1{(3-2)\cdot2}
=\frac12+\frac14+\frac12=\frac54.
\]
\end{example}

Let $\lambda/\mu$ be a skew Young diagram and $N$ a positive integer. For an index $\bk$ on $\lambda/\mu$, the \emph{Schur multiple harmonic sum} is
\begin{equation}\label{eq:schur-mhs}
\zeta_{<N}(\bk)=\sum_{T\in\SSYT_{<N}(\lambda/\mu)}\prod_{(i,j)\in\lambda/\mu}\frac{1}{T(i,j)^{k_{i,j}}}.
\end{equation}

For integers $a,b$, put $[a,b]_\Z=\{n\in\Z\mid a\leq n\leq b\}$. Following Yamamoto \cite[Definition~3.1]{Y2}, a pair $(J,J')$ of intervals in $\Z$ is called \emph{consecutive} if $J=[j_0,j_1]_\Z$ and $J'$ is one of $J$, $J\setminus\{j_1\}$, $J\cup\{j_0-1\}$ and $(J\cup\{j_0-1\})\setminus\{j_1\}$. For empty $J$, this means that $J'$ is empty or a singleton. For a skew Young diagram, the pairs $(J_{g},J_{{g}+1})$ of column sets in \eqref{eq:diagonal} are consecutive. For a tuple $\bm=(m_j)_{j\in J}$ of real numbers, we write
\[
\Pi(\bm)=\prod_{j\in J}m_j.
\]
For a consecutive pair $(J,J')$ and tuples $\bm=(m_j)_{j\in J}\in\R^J$ and $\bn=(n_j)_{j\in J'}\in\R^{J'}$, we define
\begin{align}\label{eq:interlacing}
\bm\vartriangleleft\bn
&\Longleftrightarrow
\begin{cases}
 m_j<n_j & \text{if }j\in J\text{ and }j\in J',\\
 n_{j-1}\leq m_j & \text{if }j\in J\text{ and }j-1\in J',
\end{cases}\\
\bm\trianglelefteq\bn
&\Longleftrightarrow
\begin{cases}
 m_j\leq n_j & \text{if }j\in J\text{ and }j\in J',\\
 n_{j-1}<m_j & \text{if }j\in J\text{ and }j-1\in J'.
\end{cases}\nonumber
\end{align}
These are the relations in \cite[Definition~3.2]{Y2}. For example, let $J=\{1,2\}$ and $J'=\{0,1\}$. For $\bm=(m_1,m_2)=(2,5)$ and $\bn=(n_0,n_1)=(1,4)$, the relation $\bm\vartriangleleft\bn$ says $n_0\leq m_1<n_1\leq m_2$, and indeed $1\leq2<4\leq5$. For $J=J'=\{1,2\}$, the relation $\bm\trianglelefteq\bn$ says $m_1\leq n_1<m_2\leq n_2$. For example, $(2,5)\trianglelefteq(4,6)$.

Suppose now that $\lambda/\mu$ is nonempty. Let $\bk$ be a diagonally constant index on $\lambda/\mu$. For $T\in\SSYT_{<N}(\lambda/\mu)$, we put $m_{{g},j}=T(j+{g},j)$ for ${g}_0\leq {g}\leq {g}_1$ and $j\in J_{g}$, and $\bm_{g}=(m_{{g},j})_{j\in J_{g}}$. Then the semi-standard Young tableau condition says exactly $\bm_{g}\vartriangleleft\bm_{{g}+1}$ for ${g}_0\leq {g}<{g}_1$. The inequalities $m_{{g},j}<m_{{g}+1,j}$ compare boxes in the same column, and the inequalities $m_{{g}+1,j-1}\leq m_{{g},j}$ compare boxes in the same row. Hence
\begin{equation}\label{eq:schur-mhs-diagonal}
\zeta_{<N}(\bk)=\sum_{\substack{\bm_{g}\in[1,N-1]_\Z^{J_{g}}\ ({g}_0\leq {g}\leq {g}_1)\\ \bm_{g}\vartriangleleft\bm_{{g}+1}\ ({g}_0\leq {g}<{g}_1)}}
\prod_{{g}={g}_0}^{{g}_1}\Pi(\bm_{g})^{-k_{g}}.
\end{equation}
Yamamoto's $\flat$-sum \cite[Definition~3.5]{Y2} expands each diagonal of exponent $k_{g}$ into $k_{g}$ layers. For $\bn=(n_j)_{j\in J}$ we write $\bN-\bn=(N-n_j)_{j\in J}$. The \emph{$\flat$-sum} is
\begin{equation}\label{eq:schur-flat}
\zeta^\flat_{<N}(\bk)
=\sum_{\substack{
\bn^{(l)}_{g}\in[1,N-1]_\Z^{J_{g}}\ ({g}_0\leq {g}\leq {g}_1,\ 1\leq l\leq k_{g})\\
\bn^{(l)}_{g}\trianglelefteq\bn^{(l+1)}_{g}\ (1\leq l<k_{g})\\
\bn^{(k_{g})}_{g}\vartriangleleft\bn^{(1)}_{{g}+1}\ ({g}_0\leq {g}<{g}_1)}}
\prod_{{g}={g}_0}^{{g}_1}
\frac{1}{\Pi(\bN-\bn^{(1)}_{g})\,\Pi(\bn^{(2)}_{g})\cdots\Pi(\bn^{(k_{g})}_{g})}.
\end{equation}

\begin{theorem}[{\cite[Theorem~3.6]{Y2}}]\label{thm:y2}
For every diagonally constant index $\bk$ and every positive integer $N$, we have
\[
\zeta_{<N}(\bk)=\zeta^\flat_{<N}(\bk).
\]
\end{theorem}
\section{A cyclic analogue of Schur multiple zeta values}\label{sec:cyclic}

In this section we define cyclic Young diagrams and cyclic Schur multiple harmonic sums, and prove \cref{thmA:finite}.

\subsection{Cyclic Young diagrams and tableaux}\label{sec:cyclic-def}

\begin{definition}\label{def:cyclic-diagram}
Let $\lambda/\mu$ be a nonempty skew Young diagram with $\mu_1<\lambda_1$ and $\mu_r=0$, where $r$ is the number of parts of $\lambda$. Let $q$ be an integer with \( q\geq\max\{\lambda'_1-\lambda'_{\lambda_1},\,\mu'_1-\mu'_{\lambda_1}\}. \) Let $\tau_q\colon\Z^2\to\Z^2$ be the translation
\[
\tau_q(i,j)=(i+q,\,j-\lambda_1).
\]
For $t\in\Z$, its $t$-th power is $\tau_q^t(i,j)=(i+tq,\,j-t\lambda_1)$. The \emph{periodic diagram} of $\lambda/\mu$ is
\[
\widetilde{(\lambda/\mu)_q}=\bigsqcup_{t\in\Z}\tau_q^t(\lambda/\mu).
\]
The \emph{cyclic Young diagram} is the quotient
\[
(\lambda/\mu)_q=\widetilde{(\lambda/\mu)_q}/\langle\tau_q\rangle.
\]
For a box $(i,j)$ of the periodic diagram, $[i,j]$ denotes its class in $(\lambda/\mu)_q$. The class $[i,j]$ is a \emph{cyclic corner} if neither the right nor the lower neighbor of $(i,j)$ belongs to the periodic diagram. Put $s=\lambda_1+q$. For $a\in\Z/s\Z$, the $a$-th \emph{cyclic diagonal class} consists of the classes $[i,j]$ with $i-j\equiv a\pmod s$. A cyclic diagonal class is \emph{empty} if it contains no box.
\end{definition}

\begin{unnumberedremark}
The translation $\tau_q$ moves a box $q$ rows down and $\lambda_1$ columns to the left. The translate $\tau_q^t(\lambda/\mu)$ lies in the columns $1-t\lambda_1,\ldots,\lambda_1-t\lambda_1$, so the union in \cref{def:cyclic-diagram} is disjoint. The quotient $(\lambda/\mu)_q$ is the set of orbits of the group $\langle\tau_q\rangle$. Each orbit contains exactly one box of $\lambda/\mu$. Hence a $\tau_q$-periodic map on the periodic diagram is determined by its values on $\lambda/\mu$.

The periodic column boundaries are $\lambda'_{j-t\lambda_1}=\lambda'_j+tq$ and $\mu'_{j-t\lambda_1}=\mu'_j+tq$ for $1\leq j\leq\lambda_1$ and $t\in\Z$. For every $j\in\Z$, column $j$ of the periodic diagram consists of the boxes $(i,j)$ with $\mu'_j<i\leq\lambda'_j$. The condition on $q$ says exactly that these boundaries are nonincreasing in $j$. Hence the periodic diagram is again a skew diagram. {The periodic diagram need not be edge-connected.} The conditions $\mu_1<\lambda_1$ and $\mu_r=0$ normalize the representative $\lambda/\mu$. Its first row is nonempty, and its last row begins in the first column. The translation $\tau_q$ changes $i-j$ by $s$, so the residue of $i-j$ modulo $s$ is well defined on $(\lambda/\mu)_q$.
\end{unnumberedremark}

\begin{definition}\label{def:cyclic-tableau}
A \emph{semi-standard cyclic Young tableau} of shape $(\lambda/\mu)_q$ is a $\tau_q$-periodic map $T\colon\widetilde{(\lambda/\mu)_q}\to\N$ with
\[
T(i,j)\leq T(i,j+1),\qquad T(i,j)<T(i+1,j)
\]
whenever both boxes belong to the periodic diagram. The set of these tableaux is denoted by $\SSYT((\lambda/\mu)_q)$, and its subset with entries less than $N$ by $\SSYT_{<N}((\lambda/\mu)_q)$. An \emph{index} on $(\lambda/\mu)_q$ is a $\tau_q$-periodic map $\bk\colon\widetilde{(\lambda/\mu)_q}\to\N$, represented by $\bk=(k_{i,j})_{(i,j)\in\lambda/\mu}$. It is \emph{diagonally constant} if $k_{i,j}$ depends only on $i-j$ modulo $s$. For $N\in\N$, its \emph{cyclic Schur multiple harmonic sum} is
\begin{equation}\label{eq:cyclic-tableau-sum}
\zeta_{<N}(\bk)=\sum_{T\in\SSYT_{<N}((\lambda/\mu)_q)}\prod_{(i,j)\in\lambda/\mu}T(i,j)^{-k_{i,j}}.
\end{equation}
\end{definition}

\begin{definition}\label{def:cyclic-flat}
Let $\bk$ be a diagonally constant index on $(\lambda/\mu)_q$. Expand each box $(i,j)$ into $k_{i,j}$ layers $(i,j,l)$ with $1\leq l\leq k_{i,j}$, as in Yamamoto's $\flat$-sum \cite[Definition~3.5]{Y2} (see also \eqref{eq:schur-flat}). These layers are identified periodically under $\tau_q$. Let $\mathcal F_{<N}(\bk)$ consist of periodic fillings $n=(n_{i,j,l})$ of the layers with $n_{i,j,l}\in\{1,\ldots,N-1\}$ satisfying
\[
\begin{aligned}
n_{i,j,l}&\leq n_{i,j,l+1}, &
 n_{i-1,j-1,l+1}&<n_{i,j,l},\\
n_{i,j,k_{i,j}}&<n_{i+1,j,1}, &
 n_{i,j-1,1}&\leq n_{i,j,k_{i,j}}
\end{aligned}
\]
whenever the indicated layers exist. The \emph{cyclic Schur $\flat$-sum} is
\begin{equation}\label{eq:cyclic-geometric-flat}
\zeta^\flat_{<N}(\bk)=\sum_{n\in\mathcal F_{<N}(\bk)}
\prod_{(i,j)\in\lambda/\mu}
\frac{1}{(N-n_{i,j,1})\,n_{i,j,2}\cdots n_{i,j,k_{i,j}}}.
\end{equation}
\end{definition}
For $k_{i,j}=1$ the product after $N-n_{i,j,1}$ in \eqref{eq:cyclic-geometric-flat} is empty.

\begin{example}
The periodicity can add a comparison or force an equality. Let $\lambda=(2,2)$ and $\mu=\varnothing$. For a semi-standard cyclic Young tableau $T$, write $m_{ij}=T(i,j)$ for $(i,j)\in\lambda/\mu$. Since $T$ is $\tau_q$-periodic, each box of the periodic diagram carries the entry of its class. In the pictures below, the shaded boxes form $\lambda/\mu$. The numbers outside the boxes are the row and column indices.

For $q=1$ we have $\tau_1(i,j)=(i+1,j-2)$. Let $a,b\geq1$ and take the index of \eqref{eq:intro-2x2},
\[
\bk=\cyctabdiag{a}{b}.
\]
Its double edges glue the right side of the box $(1,2)$ to the left side of the box $(2,1)$. Rows $0$ to $3$ of the periodic diagram are
\[
\begin{resizedtableau}
\none & \none[\color{black}\scriptstyle -1] & \none[\color{black}\scriptstyle 0] & \none[\color{black}\scriptstyle 1] & \none[\color{black}\scriptstyle 2] & \none[\color{black}\scriptstyle 3] & \none[\color{black}\scriptstyle 4]\\
\none[\color{black}\scriptstyle 0] & \none & \none & \none & \none & \periodicEntry{11} & \periodicEntry{12}\\
\none[\color{black}\scriptstyle 1] & \none & \none & *(gray!25) \periodicEntry{11} & *(gray!25) \periodicEntry{12} & \periodicEntry{21} & \periodicEntry{22}\\
\none[\color{black}\scriptstyle 2] & \periodicEntry{11} & \periodicEntry{12} & *(gray!25) \periodicEntry{21} & *(gray!25) \periodicEntry{22}\\
\none[\color{black}\scriptstyle 3] & \periodicEntry{21} & \periodicEntry{22}
\end{resizedtableau}
\]
The rows give $m_{11}\leq m_{12}\leq m_{21}\leq m_{22}$, and the columns give $m_{11}<m_{21}$ and $m_{12}<m_{22}$. These are the conditions in \eqref{eq:intro-2x2}. The periodicity adds $m_{12}\leq m_{21}$.

For $q=0$ we have $\tau_0(i,j)=(i,j-2)$. Take the index
\[
\bk=\cyctabrow{a}{b}.
\]
Its double edges glue each row to itself. The periodic diagram consists of the rows $1$ and $2$:
\[
\begin{resizedtableau}
\none & \none & \none[\color{black}\scriptstyle -1] & \none[\color{black}\scriptstyle 0] & \none[\color{black}\scriptstyle 1] & \none[\color{black}\scriptstyle 2] & \none[\color{black}\scriptstyle 3] & \none[\color{black}\scriptstyle 4]\\
\none[\color{black}\scriptstyle 1] & \none[\color{black}\cdots] & \periodicEntry{11} & \periodicEntry{12} & *(gray!25) \periodicEntry{11} & *(gray!25) \periodicEntry{12} & \periodicEntry{11} & \periodicEntry{12} & \none[\color{black}\cdots]\\
\none[\color{black}\scriptstyle 2] & \none[\color{black}\cdots] & \periodicEntry{21} & \periodicEntry{22} & *(gray!25) \periodicEntry{21} & *(gray!25) \periodicEntry{22} & \periodicEntry{21} & \periodicEntry{22} & \none[\color{black}\cdots]
\end{resizedtableau}
\]
The first row gives $m_{11}\leq m_{12}\leq m_{11}$, so $m_{11}=m_{12}$. Likewise $m_{21}=m_{22}$. Hence $\zeta_{<N}(\bk)=\zeta_{<N}(a+b,a+b)$.

For a single column of height $r$ with $q=r-1$, there is no additional condition from the periodicity. The cyclic sum is the multiple harmonic sum of its column index.
\end{example}

\begin{theorem}\label{thm:finite}
For every diagonally constant index $\bk$ on a cyclic Young diagram and every positive integer $N$,
\[
\zeta_{<N}(\bk)=\zeta^\flat_{<N}(\bk).
\]
\end{theorem}
This is \cref{thmA:finite}.

\subsection{Proof of Main Theorem A}\label{sec:transfer}\label{sec:datum}\label{sec:linear-algebra}
In this subsection we introduce diagonal coordinates and prove \cref{thm:finite}. Yamamoto proved his formula \cite[Theorem~3.6]{Y2} (see also \cref{thm:y2}) by the method of connected sums. The proof below is a matrix version of his argument.

\begin{definition}\label{def:datum}
Let $s,h\in\N$ with $h\leq s$. A \emph{cyclic datum} is
\[
\sD=\bigl((J_{g},k_{g})_{1\leq {g}\leq s};h\bigr),
\]
where $J_{g}\subset\Z$ is a finite interval, $k_{g}\in\N$, and $(J_{g},J_{{g}+1})$ is consecutive for $1\leq {g}\leq s$. Here $J_{s+1}=J_1-h$. Put $k_{s+1}=k_1$.
\end{definition}

Let $\bk$ be a diagonally constant index on $(\lambda/\mu)_q$. Number the diagonals of the periodic diagram by ${g}=i-j+\lambda_1$. The ${g}$-th diagonal belongs to the $a$-th cyclic diagonal class with $a\equiv {g}-\lambda_1\pmod s$. For ${g}\in\Z$ set
\begin{equation}\label{eq:index-to-datum}
J_{g}=\{j\in\Z\mid({g}+j-\lambda_1,j)\in\widetilde{(\lambda/\mu)_q}\},
\end{equation}
and let $k_{g}$ be the common exponent on this diagonal. If $J_{g}=\varnothing$, {put $k_{g}=1$}. In terms of the periodic column boundaries,
\[
J_{g}=\{j\in\Z\mid \mu'_j-j<{g}-\lambda_1\leq\lambda'_j-j\}.
\]
Both integer sequences $\mu'_j-j$ and $\lambda'_j-j$ are strictly decreasing. Hence each $J_{g}$ is an interval and $(J_{g},J_{{g}+1})$ is consecutive. Both sequences decrease by $s$ when $j$ increases by $\lambda_1$, so each $J_{g}$ is finite. Moreover, $J_{{g}+s}=J_{g}-\lambda_1$ and $k_{{g}+s}=k_{g}$. {Since $\lambda_1\leq s$, we obtain the associated cyclic datum}
\begin{equation}\label{eq:associated-datum}
\sD=\bigl((J_{g},k_{g})_{1\leq {g}\leq s};\lambda_1\bigr).
\end{equation}

For $T\in\SSYT((\lambda/\mu)_q)$ and ${g}\in\Z$, put $\bm_{g}=(m_{{g},j})_{j\in J_{g}}$ with
\begin{equation}\label{eq:tableau-to-data}
m_{{g},j}=T({g}+j-\lambda_1,j)\qquad(j\in J_{g}).
\end{equation}
The periodicity of $T$ gives $m_{{g}+s,j-\lambda_1}=m_{{g},j}$. The semi-standard condition gives $\bm_{g}\vartriangleleft\bm_{{g}+1}$ for all ${g}\in\Z$. Conversely, every family $(\bm_{g})_{{g}\in\Z}$ of tuples $\bm_{g}\in\N^{J_{g}}$ with these two properties comes from a unique $T$.

\begin{example}\label{ex:datum-2x2}
The cyclic square $((2,2)/\varnothing)_0$ with $\bk=\cyctabrow{2}{2}$ has $s=h=2$ and
\[
J_1=\{2,3\},\qquad J_2=\{1,2\},\qquad
(k_1,k_2)=(2,2),\quad J_3=J_1-2=\{0,1\}.
\]
Hence its cyclic datum is $\sD=((J_1,2),(J_2,2);2)$. Write $\bm_1=(x,y)$ on $J_1$ and $\bm_2=(u,v)$ on $J_2$. The periodicity $m_{{g}+s,j-\lambda_1}=m_{{g},j}$ gives $m_{3,0}=m_{1,2}$ and $m_{3,1}=m_{1,3}$. Hence $\bm_3=(x,y)$, and the two interlacing conditions $\bm_1\vartriangleleft\bm_2$ and $\bm_2\vartriangleleft\bm_3$ become
\[
u\leq x<v\leq y,\qquad x\leq u<y\leq v.
\]
They force $x=u$ and $y=v$, while $x<y$ remains. Thus the two entries in each cyclic row are equal.
\end{example}

\begin{definition}\label{def:finite-sums}
Let $\sD=((J_{g},k_{g})_{1\leq {g}\leq s};h)$ be a cyclic datum and let $N$ be a positive integer. Define
\[
\zeta_{<N}(\sD)
=\sum_{\substack{\bm_{g}\in[1,N-1]_\Z^{J_{g}}\ (1\leq {g}\leq s)\\
\bm_{g}\vartriangleleft\bm_{{g}+1}\ (1\leq {g}\leq s)}}
\prod_{{g}=1}^s\Pi(\bm_{g})^{-k_{g}},
\]
and
\[
\zeta^\flat_{<N}(\sD)
=\sum_{\substack{
\bn_{g}^{(l)}\in[1,N-1]_\Z^{J_{g}}\ (1\leq {g}\leq s,\ 1\leq l\leq k_{g})\\
\bn_{g}^{(l)}\trianglelefteq\bn_{g}^{(l+1)}\ (1\leq l<k_{g})\\
\bn_{g}^{(k_{g})}\vartriangleleft\bn_{{g}+1}^{(1)}\ (1\leq {g}\leq s)}}
\prod_{{g}=1}^s
\frac{1}{\Pi(\bN-\bn_{g}^{(1)})\,\Pi(\bn_{g}^{(2)})\cdots\Pi(\bn_{g}^{(k_{g})})}.
\]
Here $\bm_{g}=(m_{{g},j})_{j\in J_{g}}$ and $\bn^{(l)}_{g}=(n^{(l)}_{{g},j})_{j\in J_{g}}$. The tuples $\bm_{s+1}$ and $\bn^{(1)}_{s+1}$ are $\bm_1$ and $\bn^{(1)}_1$ with the indices shifted by $h$,
\begin{equation}\label{eq:seam}
m_{s+1,j}=m_{1,j+h},\qquad n^{(1)}_{s+1,j}=n^{(1)}_{1,j+h}\qquad(j\in J_{s+1}).
\end{equation}
\end{definition}
If $\sD$ is the datum associated with $\bk$, then \eqref{eq:tableau-to-data} gives
\begin{equation}\label{eq:geometric-datum-sum}
\zeta_{<N}(\sD)=\zeta_{<N}(\bk),
{\qquad\zeta^\flat_{<N}(\sD)=\zeta^\flat_{<N}(\bk)}.
\end{equation}
{For $J_{g}=\varnothing$, the tuples on that diagonal are empty and contribute weight $1$. It therefore suffices to prove $\zeta_{<N}(\sD)=\zeta^\flat_{<N}(\sD)$.}

The tuples in \cref{def:finite-sums} are strictly increasing in their column coordinate $j$. Indeed, let $j,j+1\in J_{g}$. Since $(J_{g},J_{{g}+1})$ is consecutive, $j\in J_{{g}+1}$. Then $\bm_{g}\vartriangleleft\bm_{{g}+1}$ gives $m_{{g},j}<m_{{g}+1,j}\leq m_{{g},j+1}$. For a $\flat$-sum layer with $k_{g}\geq2$, an internal relation $\ba\trianglelefteq\bb$ gives $a_j\leq b_j<a_{j+1}\leq b_{j+1}$. If $k_{g}=1$, the boundary relation gives the same conclusion. Thus all tuples in these sums lie in the sets $\mathcal S_{g}$ defined below.

Fix a diagonally constant index $\bk$ on a cyclic Young diagram and let $\sD$ be its associated datum. Let $N\geq2$ and $d_{g}=|J_{g}|$. For $1\leq {g}\leq s+1$, put
\[
\mathcal S_{g}=\{(m_j)_{j\in J_{g}}\in[1,N-1]_\Z^{J_{g}}\mid m_j<m_{j'}\text{ whenever }j<j'\},
\]
and identify $\mathcal S_{s+1}$ with $\mathcal S_1$ by \eqref{eq:seam}. Rows and columns of the matrices below are indexed by these tuples in lexicographic order. {For $J_{g}=\{j_1<\cdots<j_{d_{g}}\}$, we read $(m_j)_{j\in J_{g}}$ as $(m_{j_1},\ldots,m_{j_{d_{g}}})$.}

Write $\mathbf1_E$ for the indicator of a condition $E$. For $1\leq {g}\leq s$, define the rectangular matrix
\[
(R_{g}^{\vartriangleleft})_{\bm,\bn}=\mathbf1_{\{\bm\vartriangleleft\bn\}}
\qquad(\bm\in\mathcal S_{g},\ \bn\in\mathcal S_{{g}+1}).
\]
For $1\leq {g}\leq s+1$ and $\bm,\bn\in\mathcal S_{g}$, define the three square matrices
\[
(R_{g}^{\trianglelefteq})_{\bm,\bn}=\mathbf1_{\{\bm\trianglelefteq\bn\}},\qquad
(W_{g})_{\bm,\bn}=\frac{\delta_{\bm,\bn}}{\Pi(\bm)},\qquad
(W^\flat_{g})_{\bm,\bn}=\frac{\delta_{\bm,\bn}}{\Pi(\bN-\bm)},
\]
where $\delta_{\bm,\bn}$ is the Kronecker delta. Thus $R_{g}^{\vartriangleleft}$ has size $\binom{N-1}{d_{g}}\times\binom{N-1}{d_{{g}+1}}$, and the other three matrices have size $\binom{N-1}{d_{g}}\times\binom{N-1}{d_{g}}$. In each product below, adjacent factors have the same intermediate index set. By \eqref{eq:seam}, the full product is a square matrix. Expanding the traces gives
\begin{align}
\zeta_{<N}(\sD)
&=\tr\Bigl(\prod_{{g}=1}^s W_{g}^{k_{g}}R_{g}^{\vartriangleleft}\Bigr),\label{eq:trace-series}\\
\zeta^\flat_{<N}(\sD)
&=\tr\Bigl(\prod_{{g}=1}^s W^\flat_{g}\bigl(R_{g}^{\trianglelefteq}W_{g}\bigr)^{k_{g}-1}R_{g}^{\vartriangleleft}\Bigr).\label{eq:trace-flat}
\end{align}
Here the products run from ${g}=1$ on the left to ${g}=s$ on the right. In \eqref{eq:trace-series}, the $(\bm,\bn)$ entry of $W_{g}^{k_{g}}R_{g}^{\vartriangleleft}$ is $\Pi(\bm)^{-k_{g}}\mathbf1_{\{\bm\vartriangleleft\bn\}}$. In \eqref{eq:trace-flat}, the factors $R_{g}^{\trianglelefteq}W_{g}$ introduce the layers $\bn_{g}^{(2)},\ldots,\bn_{g}^{(k_{g})}$, and the factor $W^\flat_{g}$ gives the weight of the first layer. The trace gives the identification \eqref{eq:seam}.

For $0\leq a,b\leq N-1$, put
\[
C_N(a,b)=\frac{\binom{b}{a}}{\binom{N-1}{a}},
\]
where $\binom ba=0$ for $a>b$. Let $(J,J')$ be a consecutive pair of intervals, and let $\bm\in[0,N-1]_\Z^J$ and $\bn\in[0,N-1]_\Z^{J'}$. Following \cite[Definition~3.10]{Y2}, define
\[
\widetilde J=J\cup J',\qquad
\widetilde m_j=\begin{cases}m_j,&j\in J,\\0,&j\notin J,\end{cases}
\qquad
\widetilde n_j=\begin{cases}n_j,&j\in J',\\N-1,&j\notin J',\end{cases}
\]
and
\begin{equation}\label{eq:general-connector}
D_N(\bm,\bn)=\det\bigl(C_N(\widetilde m_i,\widetilde n_j)\bigr)_{i,j\in\widetilde J},
\end{equation}
with both index sets in increasing order. The definition \eqref{eq:general-connector} allows zero coordinates, as needed in \cref{lem:yamamoto-transport}.

\begin{lemma}[{\cite[Lemmas~3.11 and~3.12]{Y2}}]\label{lem:yamamoto-transport}
Let $N\geq2$, let $(J,J')$ be consecutive, and let $\bm=(m_j)_{j\in J}$ and $\bn=(n_j)_{j\in J'}$ be nondecreasing tuples. If $J=\varnothing$ or $J'=\varnothing$, then $D_N(\bm,\bn)=1$. Moreover, the following identities hold.
\begin{enumerate}[(i)]
\item If $\bm\in[1,N-1]_\Z^J$ and $\bn\in[0,N-1]_\Z^{J'}$, then
\begin{equation}\label{eq:y-transport-one}
\frac{D_N(\bm,\bn)}{\Pi(\bm)}
=\sum_{\substack{\bb\in[1,N-1]_\Z^J\\\bb\trianglelefteq\bn}}
\frac{D_N(\bm,\bb)}{\Pi(\bb)}.
\end{equation}
\item If $\bm\in[0,N-1]_\Z^J$ and $\bn\in[1,N-1]_\Z^{J'}$, then
\begin{equation}\label{eq:y-transport-two}
\sum_{\substack{\ba\in[1,N-1]_\Z^{J'}\\\bm\vartriangleleft\ba}}
\frac{D_N(\ba,\bn)}{\Pi(\bn)}
=\frac{D_N(\bm,\bn-\mathbf1)}{\Pi(\bN-\bn)}.
\end{equation}
\end{enumerate}
In the sums, $\bb$ and $\ba$ run over all tuples with the indicated relation.
\end{lemma}

Let $C=(C_N(a,b))_{1\leq a,b\leq N-1}$ and set
\begin{equation}\label{eq:connector}
(\Delta_{g})_{\bm,\bn}=\det\bigl(C_N(m_{j_a},n_{j_b})\bigr)_{a,b=1}^{d_{g}}
\quad(\bm,\bn\in\mathcal S_{g}),
\end{equation}
where $J_{g}=\{j_1<\cdots<j_{d_{g}}\}$. For $J=J'=J_{g}$, \eqref{eq:general-connector} gives $(\Delta_{g})_{\bm,\bn}=D_N(\bm,\bn)$. {The matrix $\Delta_{g}$ is upper triangular. Indeed, if $\bm$ follows $\bn$ in lexicographic order and $a$ is their first differing position, then $m_{j_i}\geq m_{j_a}>n_{j_a}\geq n_{j_b}$ for $i\geq a$ and $b\leq a$. Hence the last $d_{g}-a+1$ rows of the minor vanish in its first $a$ columns. {In the determinant expansion, the last $d_{g}-a+1$ rows can use only the last $d_{g}-a$ columns, so every term vanishes.} The diagonal entry at $\bm$ is $\prod_{a=1}^{d_{g}}C_N(m_{j_a},m_{j_a})>0$. Thus $\Delta_{g}$ is invertible whenever $\mathcal S_{g}$ is nonempty.} The matrices $W_{g}$ and $W^\flat_{g}$ are diagonal with positive entries. For $d_{g}=0$, all three matrices are $(1)$.

\begin{lemma}\label{lem:transport}
For $1\leq {g}\leq s$, we have
\begin{align}
W_{g}\Delta_{g}&=\Delta_{g}W_{g}R_{g}^{\trianglelefteq},\label{eq:transport-internal}\\
\Delta_{g}W_{g}R_{g}^{\vartriangleleft}W^\flat_{{g}+1}
&=W_{g}R_{g}^{\vartriangleleft}\Delta_{{g}+1}W_{{g}+1}.\label{eq:transport-boundary}
\end{align}
\end{lemma}
\begin{proof}
{For the first identity, take $J=J'=J_{g}$ in \cref{lem:yamamoto-transport}(i). Every tuple $\bb\trianglelefteq\bn$ in that sum is strictly increasing, since $b_j\leq n_j<b_{j+1}$ for adjacent coordinates. Thus the entrywise identity is}
\[
\frac{D_N(\bm,\bn)}{\Pi(\bm)}
=\sum_{\substack{\bb\in\mathcal S_{g}\\\bb\trianglelefteq\bn}}
\frac{D_N(\bm,\bb)}{\Pi(\bb)}.
\]
By \eqref{eq:connector} and \eqref{eq:general-connector}, this is the first identity.

For the second identity, put $J=J_{g}$ and $J'=J_{{g}+1}$. Since the coordinates are integers, $\bb\vartriangleleft\bn$ is equivalent to $\bb\trianglelefteq(\bn-\mathbf1)$ by \eqref{eq:interlacing}. {Apply \cref{lem:yamamoto-transport}(i) with $\bn-\mathbf1\in[0,N-1]_\Z^{J'}$ and then (ii) with $\bn$. This gives}
\begin{align*}
\sum_{\substack{\bb\in\mathcal S_{g}\\\bb\vartriangleleft\bn}}
\frac{D_N(\bm,\bb)}{\Pi(\bb)\Pi(\bN-\bn)}=\frac{D_N(\bm,\bn-\mathbf1)}{\Pi(\bm)\Pi(\bN-\bn)}=\frac1{\Pi(\bm)}\sum_{\substack{\ba\in\mathcal S_{{g}+1}\\\bm\vartriangleleft\ba}}
\frac{D_N(\ba,\bn)}{\Pi(\bn)}.
\end{align*}
{The outer expressions are the required matrix entries.} The summation tuples lie in $\mathcal S_{g}$ and $\mathcal S_{{g}+1}$, respectively. Indeed, the comparisons imply $b_j<n_j\leq b_{j+1}$ and $a_j\leq m_{j+1}<a_{j+1}$ whenever adjacent coordinates occur. Empty intervals are covered by the empty-product convention and \cref{lem:yamamoto-transport}.
\end{proof}

\begin{proof}[Proof of \cref{thm:finite}]
For $N=1$, both sums are $1$ when all $J_{g}$ are empty and $0$ otherwise. Suppose $N\geq2$. If $d_{g}>N-1$ for some ${g}$, then $\mathcal S_{g}$ is empty and both sums vanish. It remains to treat the case $d_{g}\leq N-1$ for all ${g}$.

Multiplying \eqref{eq:transport-internal} on the right by $W_{g}$ and on the left by $(\Delta_{g}W_{g})^{-1}$ gives
\[
R_{g}^{\trianglelefteq}W_{g}=(\Delta_{g}W_{g})^{-1}W_{g}(\Delta_{g}W_{g}).
\]
Combining this with \eqref{eq:transport-boundary} gives
\begin{align*}
W^\flat_{g}\bigl(R_{g}^{\trianglelefteq}W_{g}\bigr)^{k_{g}-1}R_{g}^{\vartriangleleft}
&=W^\flat_{g}(\Delta_{g}W_{g})^{-1}W_{g}^{k_{g}-1}\Delta_{g}W_{g}R_{g}^{\vartriangleleft}\\
&=W^\flat_{g}(\Delta_{g}W_{g})^{-1}W_{g}^{k_{g}}R_{g}^{\vartriangleleft}
\Delta_{{g}+1}W_{{g}+1}(W^\flat_{{g}+1})^{-1}.
\end{align*}
In the product over ${g}=1,\ldots,s$, the adjacent factors cancel. Hence
\begin{align*}
\prod_{{g}=1}^s W^\flat_{g}\bigl(R_{g}^{\trianglelefteq}W_{g}\bigr)^{k_{g}-1}R_{g}^{\vartriangleleft}=W^\flat_1(\Delta_1W_1)^{-1}
\Bigl(\prod_{{g}=1}^s W_{g}^{k_{g}}R_{g}^{\vartriangleleft}\Bigr)
\Delta_{s+1}W_{s+1}(W^\flat_{s+1})^{-1}.
\end{align*}
Under \eqref{eq:seam}, $W_{s+1}=W_1$, $W^\flat_{s+1}=W^\flat_1$ and $\Delta_{s+1}=\Delta_1$. Put $P=W^\flat_1(\Delta_1W_1)^{-1}$. Then the right side equals
\[
P\Bigl(\prod_{{g}=1}^s W_{g}^{k_{g}}R_{g}^{\vartriangleleft}\Bigr)P^{-1}.
\]
Since {similar} matrices have the same trace, \eqref{eq:trace-series} and \eqref{eq:trace-flat} give $\zeta_{<N}(\sD)=\zeta^\flat_{<N}(\sD)$. This proves \cref{thm:finite}.
\end{proof}

\begin{example}[The matrix calculation for a cyclic $2\times2$ square]\label{ex:transfer-2x2}
Take the cyclic index $\bk=\cyctabrow{2}{2}$ on $((2,2)/\varnothing)_0$ and let $N=4$. By \cref{ex:datum-2x2}, the two entries in each row are equal, so
\begin{equation}\label{eq:example-cyclic-44}
\zeta_{<4}(\bk)=\sum_{1\leq x<y\leq3}\frac{1}{x^4y^4}
=\frac1{16}+\frac1{81}+\frac1{1296}=\frac{49}{648}.
\end{equation}
The matrices below are those in the proof of \cref{thm:finite}.

\smallskip\noindent\textit{Step 1: the matrices.} The cyclic datum is that of \cref{ex:datum-2x2}. In lexicographic order, the rows and columns are indexed by
\[
\mathcal S_1=\mathcal S_2=\{(1,2),(1,3),(2,3)\}.
\]
From the relations $\vartriangleleft$ and $\trianglelefteq$ on these tuples, we obtain
\[
R_1^{\vartriangleleft}=R_2^{\vartriangleleft}
=\begin{pmatrix}1&0&0\\1&1&0\\0&1&1\end{pmatrix},
\qquad
R_1^{\trianglelefteq}=R_2^{\trianglelefteq}
=\begin{pmatrix}1&1&0\\0&1&1\\0&0&1\end{pmatrix}.
\]
For instance, the $( (1,3),(1,2) )$ entry of $R_1^{\vartriangleleft}$ is $1$, since $(1,3)\vartriangleleft(1,2)$ reads $1\leq1<2\leq3$. {The column tuple is indexed by $J_2$ in $R_1^{\vartriangleleft}$ and by $J_1$ in $R_1^{\trianglelefteq}$.} The {corresponding} entry of $R_1^{\trianglelefteq}$ is $0$, since $(1,3)\trianglelefteq(1,2)$ would {require} $3\leq2$.

\smallskip\noindent\textit{Step 2: weights and {connector matrix}.} The diagonal weights are
\[
W_1=W_2=\operatorname{diag}\left(\frac12,\frac13,\frac16\right),
\qquad
W_1^\flat=W_2^\flat=\operatorname{diag}\left(\frac16,\frac13,\frac12\right).
\]
Using \eqref{eq:connector} with $N=4$ {gives $\Delta_1=\Delta_2$ from the $2\times2$ minors of $C$}
\[
C=\begin{pmatrix}
\frac13&\frac23&1\\0&\frac13&1\\0&0&1
\end{pmatrix},\qquad
\Delta_1=\Delta_2
=\begin{pmatrix}
\frac19&\frac13&\frac13\\
0&\frac13&\frac23\\
0&0&\frac13
\end{pmatrix}.
\]
For example,
\[
(\Delta_1)_{(1,2),(1,3)}
=\det\begin{pmatrix}\frac13&1\\0&1\end{pmatrix}=\frac13.
\]

\smallskip\noindent\textit{Step 3: the two transport identities.} Here \eqref{eq:transport-internal} and \eqref{eq:transport-boundary} read
\[
W_1\Delta_1=\Delta_1W_1R_1^{\trianglelefteq}
=\begin{pmatrix}
\frac1{18}&\frac16&\frac16\\
0&\frac19&\frac29\\
0&0&\frac1{18}
\end{pmatrix}
\]
and
\[
\Delta_1W_1R_1^{\vartriangleleft}W_2^\flat
=W_1R_1^{\vartriangleleft}\Delta_2W_2
=\begin{pmatrix}
\frac1{36}&\frac1{18}&\frac1{36}\\
\frac1{54}&\frac2{27}&\frac1{18}\\
0&\frac1{54}&\frac1{36}
\end{pmatrix}.
\]
Since $k_1=k_2=2$, the calculation in the proof of \cref{thm:finite} gives
\[
W_{g}^\flat R_{g}^{\trianglelefteq}W_{g}R_{g}^{\vartriangleleft}
=W_{g}^\flat(\Delta_{g}W_{g})^{-1}W_{g}^2R_{g}^{\vartriangleleft}
\Delta_{{g}+1}W_{{g}+1}(W_{{g}+1}^\flat)^{-1}
\]
for ${g}=1,2$. Hence the two products in the traces are {similar}.

\smallskip\noindent\textit{Step 4: the traces.} The two factors for ${g}=1$ are
\[
W_1^2R_1^{\vartriangleleft}
=\begin{pmatrix}
\frac14&0&0\\\frac19&\frac19&0\\0&\frac1{36}&\frac1{36}
\end{pmatrix},\qquad
W_1^\flat R_1^{\trianglelefteq}W_1R_1^{\vartriangleleft}
=\begin{pmatrix}
\frac5{36}&\frac1{18}&0\\
\frac19&\frac16&\frac1{18}\\
0&\frac1{12}&\frac1{12}
\end{pmatrix}.
\]
Under the identification \eqref{eq:seam}, the factors for ${g}=2$ are equal to those for ${g}=1$. Hence
\[
\begin{aligned}
\zeta_{<4}(\bk)
&=\tr\!\left((W_1^2R_1^{\vartriangleleft})^2\right)
=\frac1{16}+\frac1{81}+\frac1{1296}\\
&=\tr\!\left((W_1^\flat R_1^{\trianglelefteq}W_1
   R_1^{\vartriangleleft})^2\right)
=\frac{11}{432}+\frac{25}{648}+\frac5{432}
=\zeta^\flat_{<4}(\bk)=\frac{49}{648}.
\end{aligned}
\]
The two products have different diagonal entries but the same trace.
\end{example}

\section[The \texorpdfstring{{integral--series identity}}{integral--series identity}]{The {integral--series identity}}\label{sec:integral}

In this section we define the admissibility condition, construct the associated $2$-labeled poset and prove \cref{thmA:integral} by letting $N\to\infty$ in \cref{thm:finite}.

Let $\bk=(k_{i,j})$ be a diagonally constant index on $(\lambda/\mu)_q$, and put $s=\lambda_1+q$. Write $\Cor(\bk)$ for the set of cyclic corners of the underlying cyclic Young diagram.

\begin{lemma}\label{lem:no-corner}
If $\Cor(\bk)=\varnothing$, then $q=0$, $\mu=\varnothing$ and $\lambda=(\lambda_1^r)$. Moreover, the row sum $K:=\sum_{j=1}^{\lambda_1}k_{i,j}$ is independent of $i$, and for every positive integer $N$,
\[
\zeta_{<N}(\bk)=\zeta_{<N}(\underbrace{K,\ldots,K}_{r}).
\]
Consequently, $\zeta_{<N}(\bk)$ converges as $N\to\infty$ if and only if $K\geq2$.
\end{lemma}
\begin{proof}
Put $w=\lambda_1$ and consider the last box $(r,\lambda_r)$ of the last row. Every box of the periodic diagram is $\tau_q^t(i,j)=(i+tq,j-tw)$ with $(i,j)\in\lambda/\mu$ and $t\in\Z$. The lower neighbor $(r+1,\lambda_r)$ is not in the periodic diagram. Otherwise $(i,j)=(r+1-tq,\lambda_r+tw)$ for some $t$, but $i>r$ for $t\leq0$ and $j>w$ for $t\geq1$. Since $[r,\lambda_r]$ is not a cyclic corner, the right neighbor $(r,\lambda_r+1)$ is in the periodic diagram. Then $(i,j)=(r-tq,\lambda_r+1+tw)$ for some $t$. For $t=0$ we get $j>\lambda_r$, and for $t\geq1$ we get $j>w$. Hence $t\leq-1$. The condition $i\leq r$ gives $q=0$, and the condition $j\geq1$ gives $t=-1$ and $\lambda_r=w$. Thus $q=0$ and $\lambda=(w^r)$. The conditions on $q$ and $\mu$ in \cref{def:cyclic-diagram} then give $0\leq\mu'_1\leq\mu'_w=0$, and hence $\mu=\varnothing$.

Since $k_{i,j}$ depends only on $i-j$ modulo $s$, the rows of exponents are cyclic shifts of each other. Moreover, the rows of every cyclic tableau $T$ satisfy $T(i,1)\leq\cdots\leq T(i,w)\leq T(i,1)$. This implies
\[
\zeta_{<N}(\bk)=\sum_{0<n_1<\cdots<n_r<N}\prod_{i=1}^r n_i^{-K}
=\zeta_{<N}(\{K\}^r).
\]
By the convergence criterion for multiple zeta values in \cref{sec:background}, the right side converges as $N\to\infty$ if and only if $K\geq2$.
\end{proof}

The index $\bk$ is called \emph{admissible} if
\[
\begin{cases}
k_{i,j}\geq2 & \text{for every }[i,j]\in\Cor(\bk),\quad \Cor(\bk)\ne\varnothing,\\[2mm]
K\geq2 & \text{if }\Cor(\bk)=\varnothing.
\end{cases}
\]
In the second case, $K$ is the row sum in \cref{lem:no-corner}.

Next we use the layers in \eqref{eq:cyclic-geometric-flat} to define a directed graph. Let $\sD=((J_{g},k_{g})_{1\leq g\leq s};\lambda_1)$ be the cyclic datum associated with $\bk$. Let $V_{\bk}$ be the set of triples $(g,j,l)$ with $1\leq g\leq s$, $j\in J_{g}$ and $1\leq l\leq k_{g}$. By \eqref{eq:tableau-to-data}, the triple $(g,j,l)$ corresponds to the layer $(g+j-\lambda_1,j,l)$ in \cref{def:cyclic-flat}. The triple $(g,j,l)$ is labeled $\bullet$ if $l=1$ and $\circ$ if $l\geq2$.

Let $G_{\bk}$ be the directed graph on $V_{\bk}$ with the arrows
\[
\begin{aligned}
(g,j,l)&\to(g,j,l+1), &\qquad (g,j-1,l+1)&\to(g,j,l) &\qquad&(1\leq l<k_{g}),\\
(g,j,k_{g})&\to(g+1,j,1), & (g+1,j-1,1)&\to(g,j,k_{g}),
\end{aligned}
\]
whenever the triples involved belong to $V_{\bk}$. Here $(s+1,j,1)$ with $j\in J_{s+1}$ stands for $(1,j+\lambda_1,1)$, as in \eqref{eq:seam}. The arrows correspond to the comparisons in \cref{def:cyclic-flat} at the same positions. An arrow is called \emph{strict} if the corresponding comparison is $<$, and \emph{weak} otherwise. The graph $G_{\bk}$ is called \emph{acyclic} if it has no directed cycle. In this case, putting $u\leq v$ for vertices $u,v$ joined by a directed path from $u$ to $v$ (including $u=v$) gives a partial order on $V_{\bk}$, the \emph{reachability order} of $G_{\bk}$.

\begin{proposition}\label{prop:local-test}
If $q>0$, then $G_{\bk}$ is acyclic. In this case, $V_{\bk}$ with the reachability order of $G_{\bk}$ is an admissible $2$-labeled poset if and only if every cyclic corner has exponent at least $2$. If $q=0$, then $\Cor(\bk)=\varnothing$. Let $K$ be the row sum in \cref{lem:no-corner}. Then $G_{\bk}$ is acyclic if and only if $K>\lambda_1$. In this case, $V_{\bk}$ with the reachability order of $G_{\bk}$ is an admissible $2$-labeled poset.
\end{proposition}
\begin{proof}
Put $H=k_1+\cdots+k_s$ and
\[
\Phi(g,j,l)=\lambda_1(k_1+\cdots+k_{g-1}+l)+Hj\qquad((g,j,l)\in V_{\bk}).
\]
The same formula gives $\Phi(s+1,j,1)=\Phi(1,j+\lambda_1,1)$, so $\Phi$ is compatible with \eqref{eq:seam}. Its increments along the four types of arrows are $\lambda_1$, $H-\lambda_1$, $\lambda_1$ and $H-\lambda_1$, in the order of the definition of $G_{\bk}$.

Suppose first that $q>0$. Then $H\geq s>\lambda_1$. Hence $\Phi$ increases along every arrow, and $G_{\bk}$ is acyclic. Every white vertex has a predecessor. A black vertex $(g,j,1)$ has no successor exactly when $[g+j-\lambda_1,j]$ is a cyclic corner of exponent $1$. Indeed, for larger exponents the next layer is a successor. For exponent $1$ the successors correspond to the right and lower neighbors. This proves the assertions for $q>0$.

Suppose now that $q=0$. Then the conditions in \cref{def:cyclic-diagram} force a rectangle without cyclic corners, and \cref{lem:no-corner} gives $H=K\geq\lambda_1$. Thus $\Phi$ proves acyclicity unless $K=\lambda_1$. In that case all exponents are $1$, and the arrows of the fourth type form a directed cycle around each row. For $K>\lambda_1$, every black vertex has a successor, since there is no cyclic corner. Hence $V_{\bk}$ is admissible.
\end{proof}

For an admissible $\bk$, define its associated $2$-labeled poset $X_{\bk}$ as follows. If $G_{\bk}$ is acyclic, let $X_{\bk}$ be $V_{\bk}$ with the reachability order of $G_{\bk}$. By \cref{prop:local-test}, the only remaining admissible case is $\Cor(\bk)=\varnothing$ and $K=\lambda_1\geq2$. In this exceptional case, let $X_{\bk}$ be the total order with label word $(\bullet\circ^{K-1})^r$.

\begin{definition}\label{def:series-integral}
Let $\bk=(k_{i,j})$ be admissible on $(\lambda/\mu)_q$. Its \emph{cyclic Schur multiple zeta series} is
\begin{equation}\label{eq:cyclic-infinite-sum}
\zeta(\bk)=\sum_{T\in\SSYT((\lambda/\mu)_q)}
\prod_{(i,j)\in\lambda/\mu}T(i,j)^{-k_{i,j}}\in[0,\infty].
\end{equation}
\end{definition}
Since the truncated tableau sets increase with $N$ and exhaust $\SSYT((\lambda/\mu)_q)$,
\begin{equation}\label{eq:series-limit}
\zeta(\bk)=\lim_{N\to\infty}\zeta_{<N}(\bk)
\end{equation}
in $[0,\infty]$.

\begin{example}[Admissible and nonadmissible cyclic indices]\label{ex:admissible-gallery}
In the first two square diagrams below the period parameter is $q=1$, and in the next two it is $q=0$:
\[
\begin{array}{c@{\qquad}c@{\qquad}c@{\qquad}c}
\cyctabdiag{2}{1}&\cyctabdiag{1}{2}&
\cyctabrow{1}{1}&\cycsingle\\
\text{admissible}&\text{not admissible}&
\text{admissible}&\text{not admissible}
\end{array}
\]
The third index has no cyclic corner. Its row sum is $K=2$, although each of its entries is $1$. The fourth index has $K=1$, and its series diverges.

A less obvious example is the skew shape $(3,3,2)/(1,0,0)$ with $q=1$. Choose the diagonally constant index
\[
\bk=\cycgallery.
\]
As an ordinary skew diagram it has two corners, $(2,3)$ and $(3,2)$. However, the periodic translate $(3,1)=\tau_1(2,4)$ shows that $(2,4)$ is a box of the periodic diagram. Thus $(2,3)$ has a right neighbor by periodicity, and it is \emph{not} a cyclic corner. Only $(3,2)$ remains a cyclic corner. Its exponent is $2$. Hence the index is admissible, although the corner $(2,3)$ has exponent $1$.
\end{example}

Next we recall two estimates of Seki and use them to prove \cref{thm:integral}.

\begin{lemma}[{\cite[Lemma~2.2]{Se}}]\label{lem:seki-estimate}
Let $r\geq1$, $a_i,b_i\in\Z_{\geq0}$, $a_1\geq1$, and $a_i+b_i\geq1$ for $1\leq i\leq r$. Put
\[
R_{<N}(\ba;\bb)=\sum_{0<n_1<\cdots<n_r<N}
\prod_{i=1}^r(N-n_i)^{-a_i}n_i^{-b_i}.
\]
As $N\to\infty$ with $\ba,\bb$ fixed, the following estimates hold.
\begin{enumerate}[(i)]
\item One has $R_{<N}(\ba;\bb)=O((\log N)^r)$.
\item If some $i\in[1,r]_\Z$ satisfies $b_i\geq1$ and $a_i+b_i\geq2$, then $R_{<N}(\ba;\bb)=O(N^{-1}(\log N)^r)$.
\item If some $i,j\in[1,r]_\Z$ satisfy $i<j$, $a_i\geq2$ and $b_j\geq1$, then $R_{<N}(\ba;\bb)=O(N^{-1}(\log N)^r)$.
\end{enumerate}
\end{lemma}

\begin{lemma}[{\cite[Proposition~2.3]{Se}}]\label{lem:seki-strict}
Let $\bk=(k_1,\ldots,k_r)$ be any nonempty index, $w=\wt(\bk)$, and let $(a_i)_{i=1}^w$ be the sequence of labels in \cref{thm:iterated-integral}. Then, without an admissibility assumption,
\[
\zeta^\flat_{<N}(\bk)
-\sum_{0<n_1<\cdots<n_w<N}\prod_{i=1}^w
\frac{1}{(N-n_i)^{a_i}n_i^{1-a_i}}
=O\bigl(N^{-1}(\log N)^w\bigr).
\]
The implied constant depends only on $\bk$.
\end{lemma}
The explicit exponent $w$ in \cref{lem:seki-strict} follows from the proof of \cite[Proposition~2.3]{Se}. In that proof the difference is a finite sum of values $R_{<N}(\ba;\bb)$ with at most $w$ summation variables. Each of them satisfies the hypothesis of \cref{lem:seki-estimate}(ii).

\begin{lemma}\label{lem:approx}
Let $X$ be an admissible $2$-labeled poset with $w$ elements. Choose a set of comparisons $n_u<n_v$ or $n_u\leq n_v$, with $u<_Xv$, which generates the order of $X$. For each positive integer $N$, put
\[
T_N=\sum_{(n_v)_{v\in X}\in[1,N-1]_\Z^X}\ \prod_{\substack{v\in X\\ \delta_X(v)=\bullet}}\frac{1}{N-n_v}\prod_{\substack{v\in X\\ \delta_X(v)=\circ}}\frac{1}{n_v},
\]
where the tuples $(n_v)_{v\in X}$ satisfy the chosen comparisons. Then $T_N\to I(X)$ as $N\to\infty$.
\end{lemma}
\begin{proof}
The case $w=0$ follows from our conventions. Suppose $w>0$, and split $T_N$ according to whether the coordinates $n_v$ are all distinct. For a linear extension $L$ of $X$, let $\bk_L$ be the admissible index given by its label word. The part with distinct coordinates is the sum over $L$ of the sums over $0<n_1<\cdots<n_w<N$ in \cref{lem:seki-strict} with $\bk=\bk_L$. By Seki's estimate \cref{lem:seki-strict}, each of them differs from $\zeta^\flat_{<N}(\bk_L)$ by a quantity tending to zero. Hence its limit is $\zeta(\bk_L)$ by the Maesaka--Seki--Watanabe formula \cref{thm:msw}.

For a tuple $(n_v)_{v\in X}$ in the remaining part, let $\rho:X\twoheadrightarrow[1,r]_\Z$ with $r<w$ map $v$ to the rank of $n_v$ among the distinct values of the coordinates. Then
\[
\begin{aligned}
\rho(u)&<\rho(v)&&\text{for every chosen comparison }n_u<n_v,\\
\rho(u)&\leq\rho(v)&&\text{for every chosen comparison }n_u\leq n_v.
\end{aligned}
\]
Conversely, each such $\rho$ and each $0<n_1<\cdots<n_r<N$ give one allowed tuple by $n_v=n_{\rho(v)}$. Let $B_h=\rho^{-1}(h)$ be the fibers of $\rho$. Put $a_h=\#\{v\in B_h:\delta_X(v)=\bullet\}$ and $b_h=\#\{v\in B_h:\delta_X(v)=\circ\}$. Each such map $\rho$ contributes the sum
\[
R_{<N}(\ba;\bb)=\sum_{0<n_1<\cdots<n_r<N}\prod_{h=1}^r\frac{1}{(N-n_h)^{a_h}n_h^{b_h}}
\]
in the notation of \cref{lem:seki-estimate}. For $v\in B_1$, choose $u\in\min X$ with $u\leq_Xv$. Then $1\leq\rho(u)\leq\rho(v)=1$, so $u\in B_1$. Likewise $B_r\cap\max X\ne\varnothing$. As $X$ is admissible, this gives $a_1\geq1$ and $b_r\geq1$. Since $r<w$, we have $\max_h(a_h+b_h)\geq2$. Choose a block $B_h$ with at least two elements. If $B_h$ contains a $\circ$-vertex, then part~(ii) of \cref{lem:seki-estimate} applies. Otherwise $a_h\geq2$ and $b_h=0$. Since $b_r\geq1$, we have $h<r$, so \cref{lem:seki-estimate}(iii) applies. In either case the contribution is $O(N^{-1}\log^rN)$ and tends to zero. There are only finitely many such maps $\rho$. Consequently, by \cref{lem:linear-extension},
\[
\lim_{N\to\infty}T_N=\sum_L\zeta(\bk_L)=I(X).\qedhere
\]
\end{proof}

\begin{theorem}\label{thm:integral}
Let $\bk$ be an admissible diagonally constant index on a cyclic Young diagram. Then the series $\zeta(\bk)$ converges{, and we have the integral--series identity}
\[
\zeta(\bk)=I(X_{\bk}).
\]
\end{theorem}
\begin{proof}
If $G_{\bk}$ is acyclic, then $X_{\bk}$ is an admissible $2$-labeled poset by \cref{prop:local-test} and the admissibility of $\bk$. Choose the comparisons $n_u<n_v$ for the strict arrows $u\to v$ of $G_{\bk}$ and $n_u\leq n_v$ for the weak ones. They generate the order of $X_{\bk}$, and by \cref{def:cyclic-flat} the sum $T_N$ of \cref{lem:approx} is then exactly $\zeta^\flat_{<N}(\bk)$. Hence \cref{lem:approx} gives
\[
I(X_{\bk})=\lim_{N\to\infty}\zeta^\flat_{<N}(\bk).
\]
By \cref{thm:finite} and \eqref{eq:series-limit},
\[
I(X_{\bk})=\lim_{N\to\infty}\zeta_{<N}(\bk)=\zeta(\bk).
\]

It remains to consider the exceptional case in which $G_{\bk}$ is not acyclic. By \cref{prop:local-test} and the admissibility of $\bk$, we have $\Cor(\bk)=\varnothing$ and $K=\lambda_1\geq2$. Then \cref{lem:no-corner} gives
\[
\zeta(\bk)=\zeta(\underbrace{K,\ldots,K}_{r}).
\]
By the definition of $X_{\bk}$ and the {classical} iterated integral expression {for multiple zeta values} \cite{LM} (see also \cref{thm:iterated-integral}), the right side is $I(X_{\bk})$.
\end{proof}
This is \cref{thmA:integral}.

\begin{example}[A six-box cyclic integral]\label{ex:integral-2x2}
By \cref{thm:integral}, the all-$2$ index on $((2,2,2)/\varnothing)_1$ satisfies
\[
\zeta\tableauparens{\cycsix}
=I\!\left(\vcenter{\hbox{\begin{xy}
{(-14,-10) \ar @{{*}-o} (-10,-6)},
{(-6,-10) \ar @{{*}-o} (-10,-6)},
{(-6,-10) \ar @{{*}-o} (-2,-6)},
{(-10,-6) \ar @{o-{*}} (-6,-2)},
{(-2,-6) \ar @{o-{*}} (-6,-2)},
{(-2,-6) \ar @{o-{*}} (2,-2)},
{(-6,-2) \ar @{{*}-o} (-2,2)},
{(2,-2) \ar @{{*}-o} (-2,2)},
{(2,-2) \ar @{{*}-o} (6,2)},
{(10,6) \ar @{{*}-o} (6,2)},
{(2,6) \ar @{{*}-o} (6,2)},
{(2,6) \ar @{{*}-o} (-2,2)},
{(14,10) \ar @{o-{*}} (10,6)},
{(6,10) \ar @{o-{*}} (10,6)},
{(6,10) \ar @{o-{*}} (2,6)},
{(14,10) \ar @{o-{*}} (-6,-2)},
{(-10,-6) \ar @{o-{*}} (10,6)},
{(6,2) \ar @{o-{*}} (-14,-10)}
\end{xy}}}\right).
\]
Each box contributes a black first layer and a white second layer. The edges encode the comparisons that define $G_{\bk}$. The three long edges come from the identification $(i,j,l)\sim(i+1,j-2,l)$.
\end{example}

\begin{remark}[Cyclic multiple zeta values]\label{rem:cyclic-hmm}
Let
\[
\bk=[(k_{1,1},\ldots,k_{1,r_1}),\ldots,
(k_{d,1},\ldots,k_{d,r_d})]
\]
be an admissible multi-index as in \cref{sec:background}, and put $r=\sum_i r_i$. Arrange $r$ boxes in $d$ vertical blocks of lengths $r_1,\ldots,r_d$ as follows. Block $1$ occupies column $d$. Block $i+1$ occupies the column to the left of block $i$ and begins in the row of the last box of block $i$. Take the period $(r-d,-d)$. The resulting cyclic Young diagram has width $d$ and one box in each cyclic diagonal class. Attach the exponent $k_{i,j}$ to the $j$-th box of block $i$. The cyclic corners are the last boxes of the blocks of length at least $2$. Hence this cyclic index is admissible. Column strictness gives $n_{i,1}<\cdots<n_{i,r_i}$. The row condition gives $n_{i,r_i}\geq n_{i+1,1}$ for $i<d$, and the periodicity gives $n_{d,r_d}\geq n_{1,1}$. Consequently, the series of \cref{def:series-integral} is exactly \eqref{eq:defcmzv}. The labels of the expanded vertical blocks are $\bullet\circ^{k_{i,1}-1}\cdots \bullet\circ^{k_{i,r_i}-1}$. The comparisons across blocks are the inequalities of Hirose, Murahara and Murakami \cite[Theorem~1]{HMM} (see also \cref{thm:hmm}). Thus \cref{thmA:integral} recovers their cyclic integral--series identity as a special case.
\end{remark}

\section{A finite decomposition formula}\label{sec:decomp}

In this section we prove \cref{thmA:decomp} by a sign-reversing involution on families of lattice paths. Gessel and Krattenthaler use a similar involution in the proof of \cite[Proposition~1]{GK}.

For a diagonally constant cyclic index $\bk$ on $(\lambda/\mu)_q$, again $s=\lambda_1+q$. For $a\in\Z$ we write $[a]_s=a+s\Z$ for the residue class of $a$ modulo $s$. For partitions $\mu\subseteq\nu$ and an $s$-periodic sequence $\bkappa=(k_d)_{d\in\Z}$ of positive integers, we write
\begin{equation}\label{eq:shape-notation}
\zeta_{<N}(\nu/\mu;\bkappa)=\sum_{T\in\SSYT_{<N}(\nu/\mu)}\prod_{(i,j)\in\nu/\mu}T(i,j)^{-k_{i-j}}.
\end{equation}

Assume that $\bk$ has no empty diagonal class. For $d\in\Z$, let $k_d$ be the common exponent of the cyclic diagonal class $[d]_s$. Thus
\begin{equation}\label{eq:kappa-period}
k_{d+s}=k_d.
\end{equation}
For $1\leq j\leq\lambda_1$, put $A_j=\mu'_j-j$, and $B_j=\lambda'_j-j$. Extend these sequences by $A_{j+\lambda_1}=A_j-s$ and $B_{j+\lambda_1}=B_j-s$. The periodic Young diagram gives
\begin{equation}\label{eq:alcove}
A_1>\cdots>A_{\lambda_1}>A_{\lambda_1+1}=A_1-s,\qquad
B_1>\cdots>B_{\lambda_1}>B_{\lambda_1+1}=B_1-s.
\end{equation}
In particular, the residue classes $[B_1]_s,\ldots,[B_{\lambda_1}]_s$ are distinct. Let
\begin{equation}\label{eq:candidates}
\mathcal A=\left\{\nu\text{ a partition}\ \middle|\
\begin{gathered}
\mu\subseteq\nu,\quad\nu_1\leq\lambda_1,\quad |\nu|=|\lambda|,\\
\{[\nu'_j-j]_s\mid1\leq j\leq\lambda_1\}
=\{[B_j]_s\mid1\leq j\leq\lambda_1\}
\end{gathered}\right\}.
\end{equation}

For $\nu\in\mathcal A$ put $\beta_j=\nu'_j-j$. There is a unique permutation $\pi_\nu\in\mathfrak S_{\lambda_1}$ such that
\begin{equation}\label{eq:permutation}
\beta_j\equiv B_{\pi_\nu(j)}\pmod s\qquad(1\leq j\leq\lambda_1).
\end{equation}

\begin{theorem}\label{thm:decomp}
Let $\bk$ be a diagonally constant index on $(\lambda/\mu)_q$ with no empty diagonal class. Then for every positive integer $N$,
\begin{equation}\label{eq:decomp}
\zeta_{<N}(\bk)=\sum_{\nu\in\mathcal A}\sgn(\pi_\nu)\,\zeta_{<N}(\nu/\mu;\bkappa).
\end{equation}
Every $\nu/\mu$ with $\nu\in\mathcal A$ has $|\lambda/\mu|$ boxes and the same multiset of exponents as $\bk$. More precisely, for every $a\in\Z/s\Z$,
\begin{equation}\label{eq:residue-count-theorem}
\#\{(i,j)\in\nu/\mu\mid i-j\equiv a\pmod s\}
=
\#\{(i,j)\in\lambda/\mu\mid i-j\equiv a\pmod s\},
\end{equation}
and all boxes in this residue class carry exponent $k_a$.
\end{theorem}
This is \cref{thmA:decomp}.

\begin{proof}[Proof of \cref{thm:decomp}]
Put $w=\lambda_1$. For $\nu\in\mathcal A$ and $T\in\SSYT_{<N}(\nu/\mu)$, set
\begin{equation}\label{eq:tableau-runners}
x_j(t)=A_j+\#\{i:(i,j)\in\nu/\mu,\ T(i,j)\leq t\},
\qquad {t\in[0,N-1]_\Z},\quad 1\leq j\leq w.
\end{equation}
The functions $x_j$ are lattice paths. Column strictness and row weak increase are equivalent to
\[
x_j(t)-x_j(t-1)\in\{0,1\}\quad(1\leq t<N),\qquad
x_1(t)>\cdots>x_w(t).
\]
Indeed, the boxes of $\mu$ together with those filled by values at most $t$ form a Young diagram. Conversely, the paths recover $T$ by $T(i,j)=\min\{1\leq t<N:i\leq x_j(t)+j\}$. The initial and final positions are $x_j(0)=A_j$ and $x_j(N-1)=\beta_j=\nu'_j-j$. To a step from $u-1$ to $u$ at time $t$ we attach the weight $t^{-k_u}$ and the pair $([u]_s,t)$. The product of these weights is the tableau weight. Thus the right side of \eqref{eq:decomp} is the signed weight sum of these strictly ordered path families, with sign $\sgn(\pi_\nu)$.

Suppose that two paths meet modulo $s$. Choose the last such time $\tau$ and the lexicographically first pair $i<j$ meeting then. Since $\nu\in\mathcal A$, the endpoint residues $[\beta_1]_s,\ldots,[\beta_w]_s$ are the distinct classes $[B_1]_s,\ldots,[B_w]_s$ in some order. Hence no two paths meet modulo $s$ at time $N-1$, and $\tau<N-1$. Write $x_i(\tau)-x_j(\tau)=ds$. Keep the paths through time $\tau$ and switch their tails after time $\tau$ by
\begin{equation}\label{eq:runner-involution}
\widehat x_i(t)=x_j(t)+ds,\qquad
\widehat x_j(t)=x_i(t)-ds\qquad(t>\tau).
\end{equation}
All other paths are unchanged. At $\tau$ the positions agree with the original ones, and every later increment is still $0$ or $1$. The multiset of positions modulo $s$ at each time is unchanged, so there is no later meeting and hence no later equality of positions. Paths with these increments cannot reverse their strict order without being equal at an integer time. The new family therefore remains strictly ordered.

Its endpoints $\widehat\beta_j$ satisfy $\widehat\beta_j\geq A_j$. They have the same sum and the same residue set as the original endpoints. Consequently $\widehat\nu'_j=\widehat\beta_j+j$ defines a partition $\widehat\nu\in\mathcal A$. The switch transposes the two endpoint residues and therefore reverses the sign. It preserves every pair $([u]_s,t)$, and hence the weight by \eqref{eq:kappa-period}. It also preserves the last meeting time and all positions at that time. Applying the same switch again recovers the original family. This is a weight-preserving, sign-reversing involution.

The families that are not cancelled have no meeting modulo $s$. They satisfy
\begin{equation}\label{eq:cyclic-path-order}
x_1(t)>\cdots>x_w(t)>x_1(t)-s\qquad(t\in[0,N-1]_\Z).
\end{equation}
Indeed, they are strictly ordered, and the last inequality holds for $t=0$ by \eqref{eq:alcove}. The integer gap $x_1(t)-x_w(t)$ changes by at most $1$ per step, so it cannot reach or exceed $s$ without a meeting modulo $s$. Their endpoints have the residues of $B_1,\ldots,B_w$ in the same cyclic order. {Choose $c$ so that $\beta_1=B_{1+c}$. The two decreasing lists $(\beta_j)_{j=1}^w$ and $(B_{j+c})_{j=1}^w$ lie in $(\beta_1-s,\beta_1]$ and have the same residue classes modulo $s$, so $\beta_j=B_{j+c}$ for every $j$. Hence}
\[
\sum_{j=1}^w\beta_j=\sum_{j=1}^w B_j-cs.
\]
Since $|\nu|=|\lambda|$, the two endpoint sums are equal. Hence $c=0$. Thus $\nu=\lambda$ and $\pi_\nu=\id$. The last inequality in \eqref{eq:cyclic-path-order} is exactly the row condition between columns $w$ and $w+1$ of the periodic diagram, after extending the paths by $x_{j+w}(t)=x_j(t)-s$. The remaining families are precisely the cyclic tableaux, all with positive sign. This proves \eqref{eq:decomp}, including the case $q=0$. For $N=1$ both sides vanish.

Finally, for $\nu\in\mathcal A$ write $\beta_j=B_{\pi_\nu(j)}+sr_j$. Equality of endpoint sums gives $\sum_j r_j=0$. For an integer representative $a$ of a residue class, the number of boxes of $\nu/\mu$ with $i-j\equiv a\pmod s$ is
\begin{align*}
\sum_{j=1}^w\left(
\left\lfloor\frac{\beta_j-a}{s}\right\rfloor-
\left\lfloor\frac{A_j-a}{s}\right\rfloor\right)
&=\sum_{j=1}^w\left(
\left\lfloor\frac{B_{\pi_\nu(j)}-a}{s}\right\rfloor+r_j-
\left\lfloor\frac{A_j-a}{s}\right\rfloor\right)\\
&=\sum_{j=1}^w\left(
\left\lfloor\frac{B_j-a}{s}\right\rfloor-
\left\lfloor\frac{A_j-a}{s}\right\rfloor\right).
\end{align*}
This is the corresponding count for $\lambda/\mu$. This proves \eqref{eq:residue-count-theorem} and the preservation of the multiset of exponents.
\end{proof}

\begin{example}\label{ex:decomp-2x2}
Let $a,b\in\N$ and take the index $\bk=\cyctabrow{a}{b}$ on $((2,2)/\varnothing)_0$. Here $s=2$, $A=(-1,-2)$ and $B=(1,0)$. The partitions in $\mathcal A$ and their signs are
\[
\begin{array}{c|c|c|c}
\nu&(\nu'_1-1,\nu'_2-2)&\pi_\nu&\sgn(\pi_\nu)\\\hline
(2,2)&(1,0)&\mathrm{id}&+1\\
(2,1,1)&(2,-1)&(12)&-1\\
(1,1,1,1)&(3,-2)&\mathrm{id}&+1
\end{array}.
\]
The periodic exponents are $k_{2d}=a$ and $k_{2d+1}=b$. Thus \cref{thm:decomp} gives
\begin{equation}\label{eq:decomp-example-2x2}
\zeta_{<N}(a+b,a+b)
=\zeta_{<N}\tableauparens{\cyctabrow{a}{b}}
=\zeta_{<N}\tableauparens{\begin{ytableau}a&b\\b&a\end{ytableau}}
-\zeta_{<N}\tableauparens{\begin{ytableau}a&b\\b\\a\end{ytableau}}
+\zeta_{<N}\tableauparens{\begin{ytableau}a\\b\\a\\b\end{ytableau}}.
\end{equation}
{Each term in the signed sum} has four boxes carrying two copies each of $a$ and $b$.
\end{example}

\noeqref{eq:flat-decomp,eq:motivic-word-decomp}
\begin{remark}[A decomposition of finite $\flat$-sums]\label{rem:flat-decomp}
Under the hypotheses of \cref{thm:decomp}, for every positive integer $N$,
\begin{equation}\label{eq:flat-decomp}
\zeta^\flat_{<N}(\bk)
=\sum_{\nu\in\mathcal A}\sgn(\pi_\nu)\,
\zeta^\flat_{<N}(\nu/\mu;\bkappa).
\end{equation}
Here $\zeta^\flat_{<N}(\nu/\mu;\bkappa)$ is the $\flat$-sum \eqref{eq:schur-flat} of the diagonally constant index $(i,j)\mapsto k_{i-j}$ on $\nu/\mu$. Indeed, \cref{thmA:finite} identifies the left side with $\zeta_{<N}(\bk)$. Applying \cref{thmA:decomp} and \cref{thm:y2} to each summand gives the right side.
\end{remark}

\begin{remark}\label{rem:layered-decomp}
Assume the hypotheses of \cref{thm:decomp}. Assume also that $\bk$ is admissible. For every $\nu\in\mathcal A$, assume that the Schur index on $\nu/\mu$ given by $\bkappa$ is admissible. A layered version of the involution in the proof of \cref{thm:decomp} gives
\begin{equation}\label{eq:motivic-word-decomp}
W(X_{\bk})=\sum_{\nu\in\mathcal A}\sgn(\pi_\nu)W(X_\nu),
\end{equation}
where $X_\nu$ is the $2$-labeled poset of Hirose, Murahara and Onozuka \cite[Theorem~1.2]{HMO} for this Schur index (see also \cref{thm:hmo}). Explicitly, $X_\nu$ is the set of layers $(i,j,l)$ with $(i,j)\in\nu/\mu$ and $1\leq l\leq k_{i-j}$, with the order generated by the comparisons between layers in \eqref{eq:schur-flat}. The layer $(i,j,1)$ is labeled $\bullet$ and the layers with $l\geq2$ are labeled $\circ$. Here
\[
W(X)=\sum_{({x}_1,\ldots,{x}_{w})\in\operatorname{LE}(X)}\delta_X({x}_1)\cdots\delta_X({x}_{w})
\]
is the formal sum of the label words of all linear extensions of $X$, with multiplicities, {where $w=|X|$}. To retain the strict and weak comparisons between layers, expand each box of exponent $k_d$ into $k_d$ layers, labeled as above. When white layers occur, list the $k_d$ layers of the classes $[d]_s$ for $d=1,\ldots,s$ in this cyclic order, with $l$ increasing. Cut this cyclic sequence between a black and a white layer. Then order the layers at each time: first the black layers in decreasing order, then the white layers in increasing order. Switching the tails at the last meeting preserves the number of each label at every original time. In this cancellation, keeping only the terms in which all layer values are distinct yields \eqref{eq:motivic-word-decomp}.
\end{remark}

\section{Applications}\label{sec:applications}
In this section we give three applications of the main theorems: the version of \cref{thmA:finite} for finite multiple zeta values, a proof of Zagier's identity from \cref{thmA:integral} by counting linear extensions of a poset, and a motivic version of \cref{thmA:decomp} {given in \cref{cor:motivic-decomp}.}

\subsection{Finite multiple zeta values}\label{sec:finite-mzv}
Let
\[
\sA=\Bigl(\prod_{{p}\ \mathrm{prime}}\mathbb F_{p}\Bigr)
\Big/\Bigl(\bigoplus_{{p}\ \mathrm{prime}}\mathbb F_{p}\Bigr),
\qquad \mathbb F_{p}=\mathbb Z/{p}\mathbb Z.
\]
For an index $\bk$, the finite multiple zeta value of Kaneko and Zagier is $\zeta_\sA(\bk)=(\zeta_{<{p}}(\bk)\bmod{p})_{p}$ \cite[Definition~7.1]{Ka}. The cyclic analogue $\zeta_\sA(\bk)$ is defined in the same way using cyclic Schur multiple harmonic sums. It includes Murahara's cyclic analogue of finite multiple zeta values \cite[Definition~1.10]{Mur} after reversing the order of the vertical blocks in the realization of \cref{rem:cyclic-hmm}.

Recall the directed graph $G_{\bk}$ on $V_{\bk}$ from \cref{sec:integral}. Let $\mathfrak P(\bk)$ be the set of ordered partitions $(B_1,\ldots,B_r)$ of $V_{\bk}$ into nonempty blocks such that every arrow of $G_{\bk}$ ends in the block where it starts or in a later block, and every strict arrow ends in a later block.

\begin{corollary}\label{cor:finite-A}
For every diagonally constant cyclic index $\bk$ on $(\lambda/\mu)_q$,
\begin{equation}\label{eq:finite-A}
\zeta_\sA(\bk)
=(-1)^{|\lambda/\mu|}\sum_{(B_1,\ldots,B_r)\in\mathfrak P(\bk)}
\zeta_\sA(|B_1|,\ldots,|B_r|).
\end{equation}
\end{corollary}
\begin{proof}
Take $N={p}$ in \cref{thm:finite}. Modulo ${p}$, $({p}-n)^{-1}\equiv-n^{-1}$ for $1\leq n<{p}$. There is one such factor for each box. This gives the sign in \eqref{eq:finite-A}. Group the layer variables according to their distinct values $1\leq a_1<\cdots<a_r<{p}$. The corresponding fibers $B_h$ form precisely the ordered partitions in $\mathfrak P(\bk)$, and their contribution is $\prod_h a_h^{-|B_h|}$. Summing over the values gives $\zeta_{<{p}}(|B_1|,\ldots,|B_r|)$. The resulting congruence for each prime ${p}$ gives the asserted identity in $\sA$.
\end{proof}

\begin{example}[A cyclic row of two boxes]\label{ex:finite-blocks}
Take $\bk=\cyctworow$ on $((2)/\varnothing)_0$. The periodicity forces the two tableau entries to be equal, so $\zeta_{\sA}(\bk)=\zeta_{\sA}(4)$. Write $x,z$ for the black layer vertices and $y,u$ for the white layer vertices. The arrows of $G_{\bk}$ are $x\to y$, $x\to u$, $z\to y$ and $z\to u$, and they are all weak. The ordered partitions in $\mathfrak P(\bk)$, grouped by block sizes, are
\[
\begin{array}{c|l|c}
\text{block sizes}&\text{one representative}&\text{multiplicity}\\\hline
(4)&(\{x,z,y,u\})&1\\
(1,3)&(\{x\},\{z,y,u\})&2\\
(3,1)&(\{x,z,y\},\{u\})&2\\
(2,2)&(\{x,z\},\{y,u\})&1\\
(1,2,1)&(\{x\},\{z,y\},\{u\})&4\\
(1,1,2)&(\{x\},\{z\},\{y,u\})&2\\
(2,1,1)&(\{x,z\},\{y\},\{u\})&2\\
(1,1,1,1)&(\{x\},\{z\},\{y\},\{u\})&4
\end{array}.
\]
Exchanging $x,z$ and $y,u$ gives the indicated multiplicities. Since $(-1)^{|\lambda/\mu|}=(-1)^2=1$, \cref{cor:finite-A} gives
\begin{equation}\label{eq:finite-six-boxes}
\begin{aligned}
\zeta_{\sA}\tableauparens{\cyctworow}=\zeta_{\sA}(4)
={}&\zeta_{\sA}(4)+2\zeta_{\sA}(1,3)
+2\zeta_{\sA}(3,1)+\zeta_{\sA}(2,2)\\
&+4\zeta_{\sA}(1,2,1)+2\zeta_{\sA}(1,1,2)\\
&+2\zeta_{\sA}(2,1,1)+4\zeta_{\sA}(1,1,1,1).
\end{aligned}
\end{equation}
For the block sizes $(1,2,1)$, the two choices of the first black vertex and the two choices of the last white vertex give the coefficient $4$.
\end{example}

\subsection{An elementary proof of Zagier's identity}\label{sec:zagier}
The identity below was conjectured by Zagier \cite{Za} and proved in \cite{BBBL01} and \cite[Theorem~1]{BBBL98}. \Cref{thmA:integral} gives another proof by counting linear extensions.
\begin{corollary}\label{cor:zagier}
For every integer $n\geq0$,
\begin{equation}\label{eq:zagier}
4^n\zeta(\{1,3\}^n)=\zeta(\{4\}^n).
\end{equation}
\end{corollary}
\begin{proof}
For $n=0$ both sides of \eqref{eq:zagier} are $1$. For $n\geq1$, let $\bk_n$ be the all-$2$ index on $((2^n)/\varnothing)_0$. It is admissible, since it has no cyclic corner and its row sum is $4$. The layers of $X_{\bk_n}$ form $2n$ pairs of incomparable elements, namely the two black layers and then the two white layers of each row. Every element of a pair lies below both elements of the next pair. Hence $X_{\bk_n}$ has $2^{2n}=4^n$ linear extensions, all with label word $(\bullet\bullet\circ\circ)^n$. Thus \cref{lem:no-corner}, \cref{thmA:integral} and \cref{lem:linear-extension} give
\[
\zeta(\{4\}^n)=\zeta\tableauparens{\cyctaballtwo}
=I\!\left(\vcenter{\hbox{\begin{xy}
{(-14,-10) \ar @{{*}-o} (-10,-6)},
{(-6,-10) \ar @{{*}-o} (-10,-6)},
{(-6,-10) \ar @{{*}-o} (-2,-6)},
{(-10,-6) \ar @{o-{*}} (-6,-2)},
{(-2,-6) \ar @{o-{*}} (-6,-2)},
{(-2,-6) \ar @{o-{*}} (2,-2)},
{(-6,-2) \ar @{{*}.o} (-2,2)},
{(2,-2) \ar @{{*}.o} (-2,2)},
{(2,-2) \ar @{{*}.o} (6,2)},
{(10,6) \ar @{{*}-o} (6,2)},
{(2,6) \ar @{{*}-o} (6,2)},
{(2,6) \ar @{{*}-o} (-2,2)},
{(14,10) \ar @{o-{*}} (10,6)},
{(6,10) \ar @{o-{*}} (10,6)},
{(6,10) \ar @{o-{*}} (2,6)},
{(14,10) \ar @{o-{*}} (2,6)},
{(-2,2) \ar @{o-{*}} (10,6)},
{(6,2) \ar @{o.{*}} (-6,-2)},
{(-10,-6) \ar @{o-{*}} (2,-2)},
{(-2,-6) \ar @{o-{*}} (-14,-10)},
{(-14,-10) \ar @/^4mm/ @{-}^{2n} (6,10)}
\end{xy}}}\right)
=4^n\zeta(\{1,3\}^n).\qedhere
\]
\end{proof}

\subsection{Motivic cyclic Schur multiple zeta values}\label{sec:app-motivic}
\mathtoolsset{showonlyrefs=false}
In this subsection we define motivic cyclic Schur multiple zeta values from the posets $X_{\bk}$ of \cref{thmA:integral}. All indices in this subsection are diagonally constant and admissible in the sense of \cref{sec:background,sec:integral}.

Let $\mathcal H$ be the algebra of motivic multiple zeta values of Brown \cite[Definition~2.1]{Br}, with period map $\operatorname{per}:\mathcal H\to\mathbb R$. For a $2$-labeled poset $X$, let $b(X)$ be its number of $\bullet$-vertices. Let $I^{\mathfrak m}(X)$ be the motivic version of $I(X)$ defined by the first named author \cite[Definition~2.9]{Fa}, using $dt/(t-1)$ at each $\bullet$-vertex. For a cyclic index $\bk$, we define
\begin{equation}\label{eq:def-motivic-cyclic}
\zeta^{\mathfrak m}_{\mathrm{cyc}}(\bk)
:=(-1)^{b(X_{\bk})}I^{\mathfrak m}(X_{\bk})\in\mathcal H.
\end{equation}
The sign converts $dt/(t-1)$ to the analytic convention $dt/(1-t)$, so \cref{thmA:integral} gives
\begin{equation}\label{eq:motivic-period}
\operatorname{per}\bigl(\zeta^{\mathfrak m}_{\mathrm{cyc}}(\bk)\bigr)
=I(X_{\bk})=\zeta(\bk).
\end{equation}
For the Schur index on $\nu/\mu$ given by $\bkappa$, let $X_\nu$ be as in \cref{rem:layered-decomp}. Set
\[
\zeta^{\mathfrak m}(\nu/\mu;\bkappa)
:=(-1)^{b(X_\nu)}I^{\mathfrak m}(X_\nu).
\]

\begin{corollary}\label{cor:motivic-decomp}

Let $\bk$ be an admissible diagonally constant index on $(\lambda/\mu)_q$ with no empty diagonal class. Suppose that, for every $\nu\in\mathcal A$, the Schur index on $\nu/\mu$ given by $\bkappa$ is admissible. Then
\begin{equation}\label{eq:motivic-decomp}
\zeta^{\mathfrak m}_{\mathrm{cyc}}(\bk)
=\sum_{\nu\in\mathcal A}\sgn(\pi_\nu)
\zeta^{\mathfrak m}(\nu/\mu;\bkappa).
\end{equation}
\end{corollary}
\begin{proof}

Since $\lambda\in\mathcal A$, the Schur index on $\lambda/\mu$ is admissible. This rules out the exceptional rectangle with all exponents equal to $1$, so $G_{\bk}$ is acyclic by \cref{prop:local-test}. Applying motivic iterated integrals to the label-word identity \eqref{eq:motivic-word-decomp} in \cref{rem:layered-decomp} gives the corresponding identity for $I^{\mathfrak m}$ \cite[Definition~2.9]{Fa}. Each word on either side has one $\bullet$ for each box, and $|\nu/\mu|=|\lambda/\mu|$ for every $\nu\in\mathcal A$. Multiplying by the common sign $(-1)^{|\lambda/\mu|}$ gives the asserted identity.
\end{proof}

By \eqref{eq:motivic-period} and the integral expression {for Schur multiple zeta values given by} Hirose, Murahara and Onozuka \cite[Theorem~1.2]{HMO} (see also \cref{thm:hmo}), taking periods recovers the $N\to\infty$ limit of \eqref{eq:decomp}.

\mathtoolsset{showonlyrefs=true}

\end{document}